\documentclass[alpha-refs]{wiley-article}

\pdfoutput=1
\usepackage{siunitx}
\usepackage{hhline}
\usepackage{lscape}
\usepackage{framed}
\usepackage{savefnmark}

\papertype{Original Article}

\title{Splitting models for multivariate count data}

\author[1\authfn{1}]{Pierre Fernique}
\author[2\authfn{1}]{Jean Peyhardi}
\author[3]{Jean-Baptiste Durand}

\contrib[\authfn{1}]{Equally contributing authors.}

\affil[1]{UMR AGAP, Cirad Inra, Univ. Montpellier, Montpellier, 34000, France}
\affil[2]{IGF, Univ. Montpellier, CNRS, INSERM, Montpellier, France.}
\affil[3]{Univ. Grenoble Alpes, CNRS, Inria, Grenoble INP${}^\star$, LJK, 38000 Grenoble, France
${}^\star$ Institute of Engineering Univ. Grenoble Alpes}

\corraddress{Jean Peyhardi, Universit\'e de Montpellier, 34000, France}
\corremail{jean.peyhardi@umontpellier.fr}

\fundinginfo{} 

\runningauthor{Pierre Fernique et al.}

\newtheorem{property}{Property}

\usepackage{glossaries}

\newacronym{psd}{PSD}{power series distribution}
\newacronym{mle}{MLE}{maximum likelihood estimator}
\newacronym{mme}{MME}{mmethod of moment estimation}
\newacronym{pgf}{pgf}{probability generative function}
\newacronym{pmf}{pmf}{probability mass function}
\newacronym{glm}{GLM}{generalized linear model}

\usepackage{amssymb}
\usepackage{amsmath}

\usepackage{subfig}
\usepackage[colorlinks=true]{hyperref}
\hypersetup{
    linkcolor=blue,          
    citecolor=blue,        
    filecolor=blue,      
    urlcolor=blue           
}
\usepackage{multirow}
\usepackage[utf8]{inputenc}
\DeclareUnicodeCharacter{2212}{-}

\begin{document}

\maketitle

\begin{abstract}

Considering discrete models, the univariate framework has been studied in depth compared to the multivariate one.
This paper first proposes two criteria to define a \textit{sensu stricto} multivariate discrete distribution.
It then introduces the class of splitting distributions that encompasses all usual multivariate discrete distributions (multinomial, negative multinomial, multivariate hypergeometric, multivariate negative hypergeometric, etc \ldots) and contains several new.
Many advantages derive from the compound aspect of splitting distributions. 
It simplifies the study of their characteristics, inferences, interpretations and extensions to regression models.
Moreover, splitting models can be estimated only by combining existing methods, as illustrated on three datasets with reproducible studies.

\keywords{singular distribution, convolution distribution, compound distribution, discrete multivariate regression}
\end{abstract}

\section{Introduction}

The analysis of multivariate count data is a crucial issue in numerous application settings, particularly in the fields of biology \citep{AL12}, ecology \citep{DK12} and econometrics \citep{Win13}.
Multivariate count data are defined as the number of items of different categories issued from sampling within a population, which individuals are grouped.
Denoting by $J$ this number of categories, multivariate count data analysis relies on modeling the joint distribution of the discrete random vector $\boldsymbol{Y} = \left(Y_1,\ldots,Y_J\right)$.
In genomics for instance, the data obtained from sequencing technologies are often summarized by the counts of DNA or RNA fragments within a genomic interval (e.g., RNA seq data).
The most usual models in this situation, are multinomial and Dirichlet multinomial regression to take account of some environmental covariate effects on these counts.
In this way, \cite{XCFL13} and \cite{CL13} studied the microbiome composition (whose output are $J$ bacterial taxa counts) and \cite{ZZZS17} studied the expression count of $J$ exon sets.

However, the multinomial and Dirichlet multinomial distributions are not appropriate for multivariate count data because of their support: the discrete simplex $\Delta_{n} = \left\{ \boldsymbol{y} \in \mathbb{N}^J : \sum_{j=1}^J y_j = n \right\}$.
This kind of distributions are said to be singular and will be denoted by $\mathcal{S}_{\Delta_n}(\boldsymbol{\theta})$.
The parameter $n$ being related to the support, is intentionally noted as an index of the distribution, distinguishing it from other parameters $\boldsymbol{\theta}$ used to described the \gls{pmf}.
In this case,  any component $Y_j$ is deterministic when the $J-1$ other components are known.
Note that initially, singular versions of some multivariate distributions have been defined by \citet{Pat68} and \citet{JP72a}.
But these distinctions were unheeded until now, leading to misuse of these distributions \citep{ZZZS17}. 
Therefore, a distribution will be considered as a $J$-multivariate distribution if 
\begin{enumerate}
  \item \label{criter1} the dimension of its support is equal to the number of variables (i.e., $\dim \{\textnormal{Supp}(\boldsymbol{Y})\}=J$).
\end{enumerate}
If no confusion could arise, $J$ will be omitted in the former notations.
Another problem that occurs when defining multivariate distributions is the independence relationships between components $Y_1,\ldots,Y_J$.
For instance, the multiple Poisson distribution described by \cite{PB67}, involves $J$ mutually independent variables.
Therefore, a multivariate distribution will be considered as a \textit{sensu stricto} multivariate distribution if: 
\begin{enumerate}[resume]
  \item \label{criter2} its probabilistic graphical model is connected (i.e., there is a path between every pair of variables).
\end{enumerate}
Additionally, such a distribution is considered as an extension of a given univariate distribution if: 
\begin{enumerate}[resume]
 \item \label{criter3} all the univariate marginal distributions belong to the same family (extension),
 \item \label{criter4} all the multivariate marginal distributions belong to the same family (natural extension).
\end{enumerate} 

Even if a singular distribution is not a \textit{sensu stricto} $J$-multivariate distribution, it is very powerful as soon as the parameter $n$ is considered as a random variable.
It then becomes a map between spaces of univariate and multivariate distributions.
Assuming that $n$ follows an univariate distribution $\mathcal{L}(\psi)$ (e.g., binomial, negative binomial, Poisson etc \ldots), the resulting compound distribution, denoted by $\mathcal{S}_{\Delta_N}(\boldsymbol{\theta}) \underset{N}{\wedge} \mathcal{L}(\psi)$, is called splitting distribution.
Under mild hypothesis, splitting distributions can be considered as \textit{sensu stricto} multivariate distributions.
They include all usual multivariate discrete distributions and several news. 
Many advantages derive from the compound aspect of splitting distributions. 
The interpretation is simply decomposed into two parts: the sum distribution (intensity of the distribution) and the singular distribution (repartition into the $J$ components).
The log-likelihood can also be decomposed according to these two parts and thus easily computed. 
All usual characteristics (support, \gls{pmf}, expectation, covariance and \gls{pgf}) are also easily obtained using this decomposition.
Finally, the generalization to regression models is naturally made by compounding a singular regression by an univariate regression.

This article is organized as follows.
In Section~\ref{section:splitting} notations used all along the paper are introduced.
The definition of singular distributions is used as a building block to introduce splitting distributions.
Positive and symmetric singular distributions are introduced, easing respectively the study of criteria \ref{criter1}-\ref{criter2} and \ref{criter3}-\ref{criter4} for resulting splitting distributions.
In Section~\ref{section:convolution} the subclass of symmetric convolution distributions is introduced (e.g., the generalized Dirichlet multinomial is a singular distribution but is not symmetric).
Additive and proportional convolution distributions are then introduced to simplify respectively the calculation of marginal distributions and the inference procedure.
Sections~\ref{section:multinomial} and \ref{section:dirichlet:multinomial} focus on splitting distributions obtained with the multinomial and the Dirichlet multinomial distributions since they are both positive and additive (e.g., the multivariate hypergeometric is additive convolution distribution but not positive). 
This leads us to precisely describe fifteen multivariate extensions (among which five are natural extensions) of usual univariate distributions giving their usual characteristics. 
In Section~\ref{section:regression}, the splitting distributions are extended to regression models.
In Section~\ref{section:application} a comparison of these regression models on two benchmark datasets and an application on a mango tree dataset are proposed.

\section{Splitting distributions}
\label{section:splitting}

\subparagraph{Notations}
All along the paper focus will be made only on non-negative discrete distributions (and regression models).  
For notational convenience, the term discrete will therefore be omitted.
Let $\left\vert \boldsymbol{Y} \right\vert = \sum_{j=1}^J Y_j$ denote the sum of the random vector $\boldsymbol{Y}$ and assume that $\left\vert \boldsymbol{Y} \right\vert \sim \mathcal{L}(\psi)$.
Let $\boldsymbol{y}=\left(y_1,\ldots,y_J\right)\in\mathbb{N}^J$ denote an outcome of the random vector $\boldsymbol{Y}$.
Let $P_{B}(A)$ denotes the conditional probability of $A$ given $B$.
Let $\operatorname{E}_{\vert\boldsymbol{Y}\vert}(\boldsymbol{Y})$ and $\operatorname{Cov}_{\vert\boldsymbol{Y}\vert}(\boldsymbol{Y})$ denote respectively the conditional expectation and covariance of the random vector $\boldsymbol{Y}$ given the sum $\vert\boldsymbol{Y}\vert$.
Let $\blacktriangle_{n}^J  = \left\{ \boldsymbol{y} \in \mathbb{N}^J : \vert \boldsymbol{y} \vert \leq n \right\}$ denote the discrete corner of the hypercube.
If no confusion could arise, $J$ will be omitted in the former notations.
Let $\binom{n}{\boldsymbol{y}}~=~n!/(n-\vert\boldsymbol{y}\vert)!\prod_{j=1}^Jy_j!$ denote the multinomial coefficient defined for $\boldsymbol{y}\in \blacktriangle_{n}$.
This notation replaces the usual notation $\binom{n}{\boldsymbol{y}}=n !/\prod_{j=1}^J y_j!$ which is defined only for $\boldsymbol{y}\in \Delta_{n}$. 
Let $(a)_{n}~=~\Gamma(a+n)/\Gamma(a)$ denote the Pochhammer symbol and 
$B(\boldsymbol{\alpha})~=~\prod_{j=1}^J \Gamma(\alpha_j) / \Gamma(\vert \boldsymbol{\alpha} \vert)$ the multivariate beta function.
Let 
$${}_{2}^JF_2 \{ (a,a^{\prime}); \boldsymbol{b}; (c,c^{\prime}); \boldsymbol{s}\}~=~\sum_{\boldsymbol{y}\in \mathbb{N}^J}  \frac{(a)_{\vert\boldsymbol{y}\vert}  (a^{\prime})_{\vert\boldsymbol{y}\vert}\prod_{j=1}^J(b_j)_{y_j}}{(c)_{\vert\boldsymbol{y}\vert} (c^{\prime})_{\vert\boldsymbol{y}\vert}} \prod_{j=1}^J \frac{s_j^{y_j}}{y_j!}$$
denote a multivariate hypergeometric function. Remark that $a^{\prime}=c^{\prime}$ lead to ${}_{1}^JF_1(a;\boldsymbol{b};c;\boldsymbol{s})$ the Lauricella’s type D function \citep{Lau93}. 
Moreover, if $J=1$ then it turns out to be the usual Gauss hypergeometric function ${}_2F_1(a;b;c;s)$ or the confluent hypergeometric ${}_1F_1(b;c;s)$.

\begin{table}
\begin{tabular}{|l|c|c|}
\cline{2-3}
\multicolumn{1}{l|}{} & Positive & Symmetric  \\
\hline 
Multinomial                       & $\times$ & $\times$  \\
Dirichlet multinomial             & $\times$ & $\times$  \\
Multivariate hypergeometric       &          & $\times$  \\
Generalized Dirichlet multinomial & $\times$ &           \\
\hline
\end{tabular}
\caption{\label{table:singular:distributions}Properties of four singular distributions.}
\end{table}

\subparagraph{Definitions}
The random vector $\boldsymbol{Y}$ is said to follow a splitting distribution\footnote{It is named splitting distribution since an outcome $y\in\mathbb{N}$ of the univariate distribution $\mathcal{L}(\psi)$ is split into the $J$ components.} if there exists a singular distribution $\mathcal{S}_{\Delta_{n}}\left(\boldsymbol{\theta}\right)$ and an univariate distribution $\mathcal{L}(\psi)$ such that $\boldsymbol{Y}$ follows the compound distribution $\mathcal{S}_{\Delta_{N}}\left( \boldsymbol{\theta}\right)  \underset{N}{\wedge}  \mathcal{L}\left(\boldsymbol{\psi}\right)$.
The \gls{pmf} is then given by  $P(\boldsymbol{Y} = \boldsymbol{y}) =  P(\vert\boldsymbol{Y}\vert=\vert\boldsymbol{y}\vert) P_{\vert\boldsymbol{Y}\vert=\vert\boldsymbol{y}\vert}(\boldsymbol{Y}=\boldsymbol{y})$ assuming that $\vert\boldsymbol{Y}\vert$ follows $\mathcal{L}(\psi)$ and $\boldsymbol{Y}$ given $\vert\boldsymbol{Y}\vert=n$ follows $\mathcal{S}_{\Delta_{n}}\left(\boldsymbol{\theta}\right)$. 
Moreover, a singular distribution is said to be:
\begin{itemize}
\item positive, if its support is the whole simplex, 
\item symmetric, if it is invariant under any permutation of its components.
\end{itemize}
If the singular distribution is symmetric, then it is possible to define a non-singular extension\footnote{The symmetry ensures that the choice of the last category to complete the vector, has no impact on the distribution.} having as support a subset of $\blacktriangle_{n}$ 
(if the singular distribution is positive, then the support of the non-singular version is exactly $\blacktriangle_{n}$).
Such a distribution for the random vector $\boldsymbol{Y}$, denoted by $\mathcal{S}_{\blacktriangle_{n}}\left( \boldsymbol{\theta},\gamma\right)$,
is such that $\left(\boldsymbol{Y}, n - \left\vert \boldsymbol{Y} \right\vert\right) \sim \mathcal{S}_{\Delta_{n}^{J+1}}\left( \boldsymbol{\theta},\gamma\right)$.
Remark that all univariate distributions bounded by $n$ (denoted by $\mathcal{L}_n(\theta)$) are non-singular distributions.
The variable $Y$ is said to follow a damage distribution\footnote{It is named damage distribution since an outcome $y\in\mathbb{N}$ of the distribution $\mathcal{L}(\psi)$ is damaged into a smaller value.} if there exists a bounded distribution $\mathcal{L}_n(\theta)$ and a distribution $\mathcal{L}(\psi)$ such that $Y$ follows the compound distribution $\mathcal{L}_{N}\left(\theta\right) \underset{N}{\wedge}\mathcal{L}\left(\boldsymbol{\psi}\right)$.

\subparagraph{Examples}
The multinomial, the Dirichlet multinomial (also known as the multivariate negative hypergeometric), the multivariate hypergeometric and the generalized Dirichlet multinomial distributions are four examples of singular distributions (see Table \ref{table:singular:distributions}). 
Contrarily to others, the multivariate hypergeometric distribution with parameters $n\in\mathbb{N}$ and $\boldsymbol{k}\in\mathbb{N}^J$, is not positive since its support is the intersection of the simplex $\Delta_n$ and the hyper-rectangle $\blacksquare_{\boldsymbol{k}}=\{\boldsymbol{y}\in\mathbb{N}^J:y_1\leq k_1,\ldots,y_J\leq k_J\}$. 
Contrarily to others, the generalized Dirichlet multinomial distribution is not symmetric.

\subsection{Splitting distributions as \textit{sensu stricto} multivariate extensions}

\subparagraph{Support}

\begin{figure}
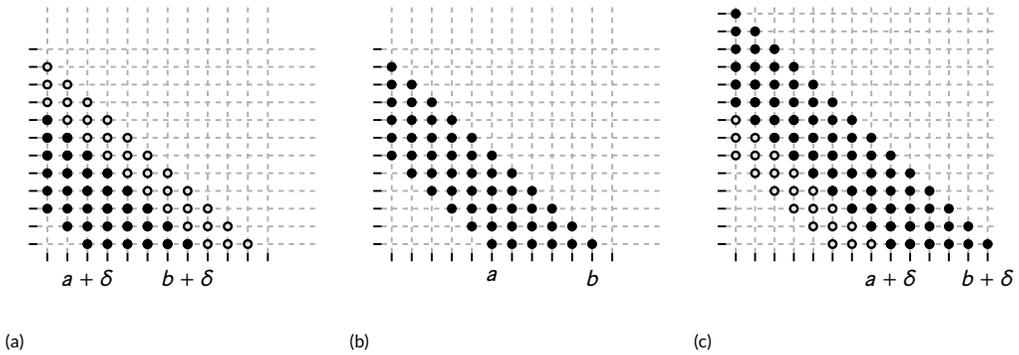

	\begin{center}
     \subfloat[][]{\scalebox{.9}{\input{support_neg.pgf}}\label{fig:splitting:support:neg}}
     \subfloat[][]{\scalebox{.9}{\input{support_null.pgf}}\label{fig:splitting:support:null}}
     \subfloat[][]{\scalebox{.9}{\input{support_pos.pgf}}\label{fig:splitting:support:pos}}
    \caption{Illustration of the effect of shifting the sum distribution on the support of positive splitting distributions.
    The support of the sum distribution is $\left\{a, \cdots, b\right\}$.
    Vectors of integers that are part of the shifted support are represented by dark circles.
    Vectors of integers that were part of the non-shifted support but are not part of the shifted support are represented by white circles.
    (a) The splitting distribution support with a negative shift (i.e., $\delta < 0$).
    (b) The splitting distribution support without any shift.
    (c) the splitting distribution support with a positive shift (i.e., $\delta > 0$).}
    \label{fig:splitting:support}
    \end{center}
\end{figure}

Firstly, let us remark that a singular distribution could be viewed as particular splitting distribution if the sum follows a Dirac distribution, i.e. $\mathcal{S}_{\Delta_{n}}\left( \boldsymbol{\theta}\right) =\mathcal{S}_{\Delta_{N}}\left( \boldsymbol{\theta}\right)  \underset{N}{\wedge}  \mathbb{1}_n$.
Assume that the dimension of a set $A\subseteq \mathbb{N}^J$ is defined as the dimension of the smaller $\mathbb{R}$-vectorial space including $A$.
The dimension of the support of a positive splitting distribution is depending on the support of the sum distribution as follows:
$$
\dim \left[ \textnormal{Supp}\left\lbrace \mathcal{S}_{\Delta_{N}}\left( \boldsymbol{\theta}\right) \underset{N}{\wedge} \mathcal{L}(\psi) \right\rbrace \right] = \left\{
\begin{array}{ll}
0    & \textnormal{if} \mathcal{L}(\psi) =  \mathbb{1}_0,\\
J-1  & \textnormal{if} \mathcal{L}(\psi) =  \mathbb{1}_n \textnormal{ with } n\in\mathbb{N}^*,\\
J    & \textnormal{otherwise}.
\end{array}
\right.
$$
Therefore all positive splitting distributions are considered as multivariate distribution (the criterion \ref{criter1} holds) when the sum is not a Dirac distribution (only non-Dirac distributions will therefore be considered all along the paper).
For non-positive splitting distributions, the dimension also depends on the support of the singular distribution.
In order to study the support of splitting distributions in a general way, 
the singular distribution is thus assumed to be positive (e.g., the case of the multivariate hypergeometric splitting distributions is omitted).
The support of $\boldsymbol{Y}$ can be expressed in terms of the sum support. 
If $\operatorname{supp}\left(\left\vert \boldsymbol{Y} \right\vert\right) = \left\{a, \cdots, b\right\}$ (with $a\in \mathbb{N}$, $b\in \mathbb{N}\cup \{\infty\}$ and $a<b$) then $\operatorname{supp}\left(\boldsymbol{Y}\right) = \blacktriangle_{b} \setminus \blacktriangle_{a - 1}$.
But, for any $j=1,\ldots,J$, the marginal support is $\operatorname{supp}\left(Y_{j}\right) = \blacktriangle_{b}$.
It can thus be interesting to consider a shifted sum $\left\vert \boldsymbol{Y} \right\vert = Z + \delta$, with $Z \sim \mathcal{L}(\boldsymbol{\psi})$, $\operatorname{supp}\left(Z\right) \subseteq  \left\{b+1, \ldots, \infty\right\}$ and $\delta \in \left\{- b, \ldots, \infty\right\} $.
As illustrated on Figure~\ref{fig:splitting:support}, such distributions can be useful for splitting distributions since they enable their supports to be modified.
In particular, for sum distributions with $a=0$ (e.g., binomial, negative binomial and Poisson distributions), $\delta = 1$ enable to remove the null vector from the joint support while it is kept in marginal supports.
Conversely, for sum distributions with $a=1$ (e.g., geometric and logarithmic series distributions), $\delta = -1$ enable to add the null vector to the joint support.
If multivariate zeros are from one structural source,
it is possible to use splitting models with shifted sums to extend hurdle models described by \citet{Mul86}.

\subparagraph{Graphical model}
A probabilistic graphical model (or graphical model, in short) is defined by a distribution and a graph such that all independence assertions that are derived from the graph using the global Markov property hold in the distribution \cite{KF13}.
A graphical model is said to be minimal, if any edge removal in the graph induces an independence assertion that is not held in the distribution.
A graphical model is said to be connected if there exists a path containing all its vertices (i.e., there is no pair of independent variables).
This is a necessary condition (criterion \ref{criter2}) to obtain a \textit{sensu stricto} multivariate distribution.
\cite{PF17} characterized the graphical model of multinomial and Dirichlet multinomial splitting distributions according to the sum distribution.
But, when the graphical model is unknown, it is sufficient to show that covariances are strictly positive to obtain a graph with at least one path between every pair of random variables.
Moments can be derived using the law of total expectation $\operatorname{E}\left(\boldsymbol{Y}\right) =  \operatorname{E}\left\lbrace \operatorname{E}_{\vert \boldsymbol{Y} \vert}\left( \boldsymbol{Y} \right) \right\rbrace$
and covariance $\operatorname{Cov}\left(\boldsymbol{Y}\right) = \operatorname{E}\left\lbrace \operatorname{Cov}_{\vert \boldsymbol{Y} \vert}\left(\boldsymbol{Y} \right)\right\rbrace  +  \operatorname{Cov}\left\lbrace \operatorname{E}_{\vert \boldsymbol{Y} \vert}\left(\boldsymbol{Y}\right)\right\rbrace$.
This method could be used for multivariate hypergeometric and generalized Dirichlet multinomial splitting distributions since their graphical models have not yet been characterized.
Similarly, the \gls{pgf} of splitting distributions can be obtain from the \gls{pgf} of the singular distribution since $G_{\boldsymbol{Y}}(\boldsymbol{s}) = \operatorname{E}\left\lbrace \bar{G}_{\boldsymbol{Y}}\left( \boldsymbol{s} \right) \right\rbrace$ where $\boldsymbol{s} = (s_1,\ldots,s_J)$ and $\bar{G}_{\boldsymbol{Y}}$ denotes the \gls{pgf} of $\boldsymbol{Y}$ given the sum $\vert \boldsymbol{Y} \vert$.

\subparagraph{Marginal distributions}
Splitting distributions that are \emph{sensu stricto} multivariate distributions (i.e., with criteria \ref{criter1} and \ref{criter2}) are not necessarily multivariate extensions.
To be considered as a multivariate extension of a specific family, marginal distributions of $Y_j$ must belong to this family.
The symmetry of the singular distribution is a sufficient condition to obtain a multivariate extension.
In fact, marginals of symmetric splitting distributions follow the same damage distribution but with different parameters (see Appendix \ref{proofs} for details).

\subsection{Inference}
If the parameters ${\boldsymbol{\theta}}$ and ${\boldsymbol{\psi}}$ are unrelated, the log-likelihood of the splitting distribution, denoted by $ \mathcal{L}\left(\boldsymbol{\theta}, \boldsymbol{\psi} ; \boldsymbol{y}\right)$, can be decomposed into log-likelihoods for the singular distribution and the sum distribution:
	\begin{align}
        \mathcal{L}\left(\boldsymbol{\theta}, \boldsymbol{\psi} ; \boldsymbol{y}\right)
        	&=  {\textstyle \log\left\{P_{\vert \boldsymbol{Y} \vert = \vert \boldsymbol{y} \vert}\left(\boldsymbol{Y} = \boldsymbol{y}  \right)\right\}} +
\log\left\{P\left(\left\vert \boldsymbol{Y} \right\vert = \vert \boldsymbol{y} \vert \right)\right\},  \nonumber \\
       \mathcal{L}\left(\boldsymbol{\theta}, \boldsymbol{\psi} ; \boldsymbol{y}\right)
            &= \mathcal{L}\left(\boldsymbol{\theta}; \boldsymbol{y} \right) +\mathcal{L}\left(\boldsymbol{\psi}; \vert \boldsymbol{y} \vert\right). \label{prop:inference:freq}
    \end{align}
Therefore, the maximum likelihood estimator (\gls{mle}) of a splitting distribution with unrelated parameters can be obtained separately using respectively the \gls{mle} of the singular distribution and the \gls{mle} of the sum distribution.
Hence, with $C$ estimators of singular distributions and $L$ estimators of univariate distributions, one is able to estimate $C \times L$ multivariate distributions, with time complexity in $\mathcal{O}\left(C + L\right)$.
Let us remark that decomposition~\eqref{prop:inference:freq} stays true for decomposable scores such as AIC and BIC.
Hence model selection using decomposable scores is also reduced to two separate model selection problems and has the same linear time and space complexity.
To limit the scope of this paper, focus will be given only to parametrization and inference of singular distributions.
Sum distributions considered in this paper are usual power series distributions and some beta compound distributions (see Appendix \ref{RDUD} for definitions and Table~\ref{table:dud} for some inference references).

\definecolor{Gray}{gray}{0.9}
\begin{table}
   \begin{tabular}{|c|c|c|c|}
     \hline
    	Distribution & Notation & Parameter Inference \\\hline\hline
    	Binomial & $\mathcal{B}_{n}\left(p\right)$ & See \cite{BD81} \\\hline
    	Negative binomial & $\mathcal{NB}\left(r, p\right)$ & See \cite{DBB12} \\\hline
        Poisson & $\mathcal{P}\left(\lambda\right)$ & See \cite{JKK93} \\\hline
    	Logarithmic series & $\mathcal{L}\left(p\right)$ & See \cite{JKK93} \\\hline\hline
     	Beta binomial & $\beta\mathcal{B}_{n}\left(a, b\right)$ & See \cite{TGRGJ94,LL12} for $n$ known \\
     	Beta negative binomial & $\beta\mathcal{NB}\left(r, a, b\right)$ & See \cite{Irw68} \\\hline
        Beta Poisson & $\beta_{\lambda}\mathcal{P}\left(a, b\right)$ & See \cite{Gur58,VWKNWRP16}\\\hline
   \end{tabular}
   \caption{\label{table:dud} References of parameter inference procedures for seven usual univariate discrete distributions.}
\end{table}

\section{Convolution splitting distributions}
\label{section:convolution}
In order to study thoroughly the graphical models and the marginals of splitting distributions, additional assumptions are necessary concerning the parametric form of the singular distribution.
Convolution splitting distributions have been introduced by \citet{S77} for $J=2$ and extended by \citet{RS79} for $J\geq 2$, but were only used as a tool for characterizing univariate discrete distributions $\mathcal{L}(\psi)$.
We here consider convolution splitting distributions as a general family of multivariate discrete distributions.

\subparagraph{Definition}
The random vector $\boldsymbol{Y}$ given $\vert \boldsymbol{Y} \vert =n$  is said to follow a (singular) convolution distribution if there exists a non-negative parametric sequence $a:=\{ a_{\theta}(y)\}_{y\in \mathbb{N}}$ such that for all $\boldsymbol{y}\in\Delta_{n}$ we have
\begin{equation*}
	P_{\vert \boldsymbol{Y} \vert =n} \left(\boldsymbol{Y} = \boldsymbol{y} \right) = \frac{1}{c_{\boldsymbol{\theta}}(n)} \prod_{j=1}^J a_{\theta_j}(y_j),
\end{equation*}
where $c_{\boldsymbol{\theta}}$ denotes the normalizing constant (i.e., the convolution of $a_{\theta_1},\ldots,a_{\theta_J}$ over the simplex $\Delta_{n}$). 
This convolution distribution, characterized by the parametric sequence $a=\{ a_{\theta}(y)\}_{y\in \mathbb{N}}$, is denoted by $\mathcal{C}_{\Delta_{n}}\left(a;\boldsymbol{\theta}\right)$ where $\boldsymbol{\theta}=(\theta_1,\ldots,\theta_J)\in\Theta^J$.
If the distribution of $\left\vert \boldsymbol{Y} \right\vert$ belongs to a family of univariate parametric distribution $\mathcal{L}(\boldsymbol{\psi})$ then, we note this compound distribution as follows
$
	\boldsymbol{Y} \sim \mathcal{C}_{\Delta_{N}}\left(a;\boldsymbol{\theta}\right)  \underset{N}{\wedge}  \mathcal{L}(\boldsymbol{\psi})
$. 
Remark that a convolution distribution is positive if and only if $a_{\theta}(y)>0$ for all $\theta \in \Theta$ and all $y \in \mathbb{N}$.
Moreover, a convolution distribution is said to be: 
\begin{itemize}
\item additive, if $c_{(\theta,\theta^{\prime})}(n) = a_{\theta + \theta^{\prime}}(n)$ for all $(\theta, \theta^{\prime}) \in \Theta^2$ and all $n \in \mathbb{N}$,
\item proportional, if $\mathcal{C}_{\Delta_{N}}\left(a;\boldsymbol{\theta}\right)=\mathcal{C}_{\Delta_{N}}\left(a;\lambda\cdot\boldsymbol{\theta}\right)$ for all $\boldsymbol{\theta}\in\Theta^J$ and all $\lambda\in\Theta$.
\end{itemize}
Thanks to the symmetry, the non-singular extension denoted by $\mathcal{C}_{\blacktriangle_{n}}\left(a;\boldsymbol{\theta},\gamma\right)$ is well defined.  Its \gls{pmf} is given by
$$
P\left(\boldsymbol{Y} = \boldsymbol{y} \right) =  \frac{1}{c_{\boldsymbol{\theta},\gamma}(n)}  a_{\gamma}(n-\vert\boldsymbol{y}\vert) \prod_{j=1}^J a_{\theta_j}(y_j),$$ for all $\boldsymbol{y}\in\blacktriangle_{n}$.
If the non-singular convolution distribution is univariate then it is denoted by $\mathcal{C}_n(a;\theta,\gamma)$.
The random variable $Y$ is said to follow a convolution damage distribution if there exists a latent variable $N$ such that $Y$ given $N=n$ follows $\mathcal{C}_{n}\left(a;\theta,\gamma\right)$ for all $n \in \mathbb{N}$.
Moreover if $N\sim\mathcal{L}\left(\boldsymbol{\psi}\right)$, we denote this compound distribution as follows:
$
    Y \sim \mathcal{C}_{N}\left(a;\theta,\gamma\right)  \underset{N}{\wedge}  \mathcal{L}\left(\boldsymbol{\psi}\right)
$.

\subparagraph{Examples}
Remark that these convolutions distributions are symmetric by construction. 
The generalized Dirichlet multinomial distribution thus not belong to this family.
The three other singular distributions are additive convolution distributions; details for the multinomial and Dirichlet multinomial distributions are given in the next two sections.
Additivity is necessary to easily determine the marginal distributions of corresponding splitting distributions.
Among those four distributions, only the multinomial distribution is a proportional convolution distribution (details are given in the next section); see Table \ref{table:convolution:distributions} for summary of these properties.
For the univariate case, the binomial, beta binomial (negative hypergeometric) and  hypergeometric distributions are three non-singular convolutions distributions.

\subsection{Properties of convolution splitting distributions}

\subparagraph{Graphical model}
Until recently, only marginal independence have been studied through the well-known Rao-Rubin condition \citep{RR64}.
\citet{PF17} generalized this condition for conditional independence and deduced the 
graphical models for positive convolution distributions under mild hypotheses.
In this case, it has been shown that there exists only one univariate distribution $\mathcal{L}(\psi)=(p_k^*)_{k\in\mathbb{N}}$ such that the minimal graphical model for $\boldsymbol{Y}$ is empty.
This distribution belongs to the family of power series distributions since there exists some $\alpha > 0$ such that $p_k^* = p_0^*  \alpha^k c_{\boldsymbol{\theta}}(k) / c_{\boldsymbol{\theta}}(0)$, for all $k \in \mathbb{N}$.
In this case and only this case, the resulting convolution splitting distribution is not a \textit{senso stricto} multivariate distribution.
For all other univariate distributions $\mathcal{L}(\psi) \neq (p_k^*)_{k\in\mathbb{N}}$, the graphical model is complete (the criterion \ref{criter2} holds) and thus the resulting convolution splitting distribution is a \emph{sensu stricto} multivariate distribution.
This characterization of the graphical model \citep[][Theorem~4]{PF17} holds for multinomial and Dirichlet multinomial splitting distributions.
But if the convolution is not positive (e.g., the multivariate hypergeometric distribution), this characterization does not hold anymore and stays an open issue. 

\subparagraph{Derived distributions}
Several stability properties hold for additive convolution splitting distributions including results of \citet{Pat68,JP72a,Xek86} as particular cases.
\begin{theorem}
\label{marginal_csd}
Let $\boldsymbol{Y}$ follow an additive convolution splitting distribution $\boldsymbol{Y}\sim\mathcal{C}_{\Delta_N}\left(a; \boldsymbol{\theta} \right)\underset{N}{\wedge}\mathcal{L}(\psi)$ then:
\begin{enumerate}
\item \label{item:marginal:sum} The marginal sum $\vert\boldsymbol{Y}_{\mathcal{I}}\vert$ follows the damage distribution $\mathcal{C}_{N}(a; \vert \boldsymbol{\theta}_{\mathcal{I}} \vert, \vert \boldsymbol{\theta}_{-\mathcal{I}} \vert ) \underset{N}{\wedge} \mathcal{L}(\boldsymbol{\psi})$.

\item \label{item:marginal:singular} The subvector $\boldsymbol{Y}_{\mathcal{I}}$ given $ \vert\boldsymbol{Y}_{\mathcal{I}}\vert=n$ follows the singular convolution distribution $\mathcal{C}_{\Delta_n}(a; \boldsymbol{\theta}_{\mathcal{I}} )$.

\item \label{item:marginal} The subvector  $\boldsymbol{Y}_{\mathcal{I}}$ follows the convolution splitting damage distribution 
$$\mathcal{C}_{\Delta_{N}}\left(a; \boldsymbol{\theta}_{\mathcal{I}}\right)  \underset{N}{\wedge}  \left\lbrace \mathcal{C}_{N^{\prime}}(a; \vert \boldsymbol{\theta}_{\mathcal{I}} \vert, \vert \boldsymbol{\theta}_{-\mathcal{I}} \vert ) \underset{N^{\prime}}{\wedge} \mathcal{L}(\boldsymbol{\psi}) \right\rbrace.$$

\item \label{item:full:conditional} The subvector $\boldsymbol{Y}_{\mathcal{I}}$ given  $\boldsymbol{Y}_{-\mathcal{I}} = \boldsymbol{y}_{-\mathcal{I}}$ follows the convolution splitting truncated and shifted distribution 
$$\mathcal{C}_{\Delta_N}(a; \boldsymbol{\theta}_{\mathcal{I}} ) \underset{N}{\wedge} \left[  TS_{\vert\boldsymbol{y}_{-\mathcal{I}}\vert}\left\lbrace \mathcal{L}(\boldsymbol{\psi}) \right\rbrace \right].$$

\item \label{item:general:conditional} The subvector $\boldsymbol{Y}_{\mathcal{I}}$ given $ \boldsymbol{Y}_{\mathcal{J}} = \boldsymbol{y}_{\mathcal{J}}$ follows the convolution splitting truncated and shifted damage  distribution 
$$\mathcal{C}_{\Delta_N}(a; \boldsymbol{\theta}_{\mathcal{I}} ) \underset{N}{\wedge} \left[ TS_{\vert\boldsymbol{y}_{\mathcal{J}}\vert} \left\lbrace \mathcal{C}_{N^{\prime}}(a; \vert \boldsymbol{\theta}_{\mathcal{I}\cup\mathcal{J}} \vert, \vert \boldsymbol{\theta}_{-\mathcal{I}\cup\mathcal{J}} \vert ) \underset{N^{\prime}}{\wedge} \mathcal{L}(\boldsymbol{\psi})\right\rbrace   \right].$$
\end{enumerate}
where
$\mathcal{I} \subset \{1,\ldots,J\}$, $-\mathcal{I}=\{1,\ldots,J\}\setminus\mathcal{I}$, $\mathcal{J}\subset-\mathcal{I}$ and $TS_{\delta}\{\mathcal{L}(\psi)\}$ denotes the truncated and shifted distribution $\mathcal{L}(\psi)$ with parameter $\delta\in\mathbb{N}$ (i.e., $X\sim TS_\delta\{\mathcal{L}(\psi)\}$ means that $P(X=x)=P_{Z\geq \delta}(Z=\delta+x)$ with $Z\sim \mathcal{L}(\psi)$).
\end{theorem}
\begin{corollary} 
\label{coro:natural:extension}
An additive convolution splitting distribution $\mathcal{C}_{\Delta_N}\left(a; \boldsymbol{\theta} \right)\underset{N}{\wedge}\mathcal{L}(\boldsymbol{\psi})$ is a natural multivariate extension of $\mathcal{L}(\boldsymbol{\psi})$ if 
the latter is stable under the damage process $\mathcal{C}_{N}(a; \vert \boldsymbol{\theta}_{\mathcal{I}} \vert, \vert \boldsymbol{\theta}_{-\mathcal{I}} \vert ) \underset{N}{\wedge} (\cdot)$ for any subset $\mathcal{I}\subset\{1,\ldots,J\}$.
Univariate marginals are thus following the distribution $\mathcal{L}(\boldsymbol{\psi}_j)$ for some $\boldsymbol{\psi}_j\in\boldsymbol{\Psi}$.
\end{corollary}
\begin{corollary} 
\label{coro:convolution:nonsingular}
The non-singular version of an additive convolution distribution is a specific convolution splitting distribution: 
    \begin{equation*}
        \mathcal{C}_{\Delta_N}\left(a;  \boldsymbol{\theta}\right) \underset{N}{\wedge} \mathcal{C}_{n}\left(a; \left\vert \boldsymbol{\theta} \right\vert, \gamma \right) = \mathcal{C}_{\blacktriangle_{n}}\left(a; \boldsymbol{\theta}, \gamma\right).
    \end{equation*}
\end{corollary}
The proof of the three previous result are given in Appendix \ref{proofs}. 
As a consequence from Corollary \ref{coro:convolution:nonsingular}, the univariate distribution $\mathcal{L}(\psi)=\mathcal{C}_{n}\left(a; \left\vert \boldsymbol{\theta} \right\vert, \gamma \right)$ is considered as the canonical case for a given convolution distribution. 
The parameters of the convolution distribution and the univariate bounded distribution  are independent given the sum $\left\vert \boldsymbol{\theta} \right\vert$.
This property can be benefited from into parameter inference procedures.
Such dependence disappears when the convolution distribution is proportional,
implying that inference procedure of the non-singular distribution is separable into two independent parts. 
Note that univariate marginals follow the convolution damage distribution: 
$
    	Y_{j} \sim  \mathcal{C}_{N}(a; \theta_{j}, \vert \boldsymbol{\theta}_{-j} \vert ) \underset{N}{\wedge} \mathcal{L}(\boldsymbol{\psi})
$. 
It is thus easy to highlight multivariate extensions studying these marginal distributions.

\begin{table}
\begin{tabular}{|l|c|c|c|}
\cline{2-4}
\multicolumn{1}{l|}{}& Positive  & Additive & Proportional \\
\hline 
Multinomial                 & $\times$ & $\times$ & $\times$  \\
Dirichlet multinomial       & $\times$ & $\times$ &           \\
Multivariate hypergeometric &          & $\times$ &           \\
\hline
\end{tabular}
\caption{\label{table:convolution:distributions}Properties differentiating three convolutions distributions.}
\end{table}

\section{Multinomial splitting distributions}
\label{section:multinomial}
In this section the multinomial distribution is introduced as a positive, additive and proportional convolution distribution.
Then, the general case of multinomial splitting distributions (i.e., for any sum distribution $\mathcal{L}(\psi)$) is addressed.
For six specific sum distributions, the usual characteristics of multinomial splitting distributions are described in Tables \ref{table:multinomial:psd} and \ref{table:multinomial:beta_compound}.
Finally, the canonical case of binomial sum distribution is detailed, with particular emphasize on parameter inference.

\subparagraph{Multinomial distribution}
Let $a_{\theta}(y)=\theta^y/y!$ be the parametric sequence that characterizes the  multinomial distribution as a convolution distribution.
It is positive since $\theta^y/y!>0$ for all $\theta\in\Theta=(0,\infty)$ and all $y\in\mathbb{N}$. 
It is additive, as a consequence from the binomial theorem: $(\theta+\gamma)^n=\sum_{y=0}^n \binom{n}{y} \theta^y \gamma^{n-y}$.
It implies, by induction on $n$, that the normalizing constant is  $c_{\boldsymbol{\theta}}(n)=a_{\vert\boldsymbol{\theta}\vert}(n)=\vert\boldsymbol{\theta}\vert^n/n!$.
The \gls{pmf} of the singular multinomial distribution is thus given by
\begin{equation}\label{pmf:singular:multinomial:convolution}
P_{\vert\boldsymbol{Y}\vert=n}\left(\boldsymbol{Y} = \boldsymbol{y}\right) = \binom{n}{\boldsymbol{y}} \prod_{j=1}^J \left(\frac{\theta_j}{\vert\boldsymbol{\theta}\vert}\right)^{y_j} \cdot \mathbb{1}_{\Delta_{n}}(\boldsymbol{y}),
\end{equation}
and is denoted by $\mathcal{M}_{\Delta_n}(\boldsymbol{\theta})$ with $\boldsymbol{\theta}\in(0,\infty)^J$.
This convolution is proportional, implying that the equivalence class of distributions $\{\mathcal{M}_{\Delta_n}(\lambda\cdot\boldsymbol{\theta}), \lambda\in(0,\infty)\}$ can be summarized by the representative element $\mathcal{M}_{\Delta_n}(\boldsymbol{\pi})$ where $\boldsymbol{\pi}=\frac{1}{\vert\boldsymbol{\theta}\vert}\cdot\boldsymbol{\theta}$.
The parameters vector $\boldsymbol{\pi}$ lies in the continuous simplex $c\Delta=\{\boldsymbol{\pi}\in(0,1)^J:\vert\boldsymbol{\pi}\vert=1\}$ and the \gls{pmf} reduces to its usual form.
The \gls{pmf} of the non-singular multinomial distribution, denoted by $\mathcal{M}_{\blacktriangle_n}(\boldsymbol{\theta},\gamma)$, is given by 
$$
P\left(\boldsymbol{Y} = \boldsymbol{y}\right) = \binom{n}{\boldsymbol{y}} \left(\frac{\gamma}{\vert\boldsymbol{\theta}\vert+\gamma}\right)^{n-\vert\boldsymbol{y}\vert} \prod_{j=1}^J \left(\frac{\theta_j}{\vert\boldsymbol{\theta}\vert+\gamma}\right)^{y_j} \cdot \mathbb{1}_{\blacktriangle_{n}}(\boldsymbol{y}).
$$
In the same way there exists a representative element $\mathcal{M}_{\blacktriangle_n}(\boldsymbol{\pi}^*,\gamma^*)$ with $(\boldsymbol{\pi}^*,\gamma^*)\in(0,1)^{J+1}$ such that $\vert\boldsymbol{\pi}^*\vert+\gamma^*=1$.
Knowing this constraint, the last parameter $\gamma^*=1-\vert\boldsymbol{\pi}^*\vert$ could be let aside to ease the notation and obtain $\mathcal{M}_{\blacktriangle_n}(\boldsymbol{\pi}^*)$ where the parameters vector $\boldsymbol{\pi}^*$ lies in the continuous corner of the open hypercube $c\blacktriangle=\{\boldsymbol{\pi}^*\in(0,1)^J:\vert\boldsymbol{\pi}^*\vert < 1\}$.
As a particular case of the non-singular multinomial distribution (when $J=1$), the binomial distribution is finally denoted by $\mathcal{B}_n(p)$ with $p\in(0,1)$ (which is also the representative element of its class).
Even if this new definition of multinomial distributions based on equivalence classes seems quite artificial, it is necessary to obtained all the properties that hold for convolution splitting distributions.
For instance Corollary~\ref{coro:convolution:nonsingular} becomes
\begin{equation}\label{msd}  
	\mathcal{M}_{\Delta_N}\left( \boldsymbol{\pi}\right) \underset{N}{\wedge} \mathcal{B}_{n}\left( p\right) = \mathcal{M}_{\blacktriangle_{n}}\left(p \cdot \boldsymbol{\pi}\right),
\end{equation}
with representative element notations (see Appendix \ref{proofs} for the proof).
A second point that is important to highlight is the difference between singular and non-singular multinomial distributions.
Contrarily to the widely held view that the multinomial distribution is the extension of the binomial distribution \citep{JKB97}, only the non-singular one  should be considered as the natural extension.
In fact, criterion \ref{criter4} does not hold for the singular multinomial distribution (multivariate marginals follow non-singular multinomial distributions).
Moreover, when confronted to multivariate counts, the usual inference of multinomial distributions \citep{JKB97,ZZZS17} is the inference of singular multinomial distributions such that $\forall n \in \mathbb{N}$ the random vector $\boldsymbol{Y}$ given $\vert \boldsymbol{Y} \vert = n$ follows $ \mathcal{M}_{\Delta_{n}}\left( \boldsymbol{\pi}\right)$.
Such a point of view therefore limits the possibility of comparing these distributions to other classical discrete multivariate distributions such as the negative multinomial distribution \citep{JKB97} or the multivariate Poisson distributions \citep{KM05} used for modeling the joint distribution of $\boldsymbol{Y}$.
The singular multinomial distribution must not be considered as a $J$-multivariate distribution since criterion \ref{criter1} does not hold.

\subparagraph{Properties of multinomial splitting distributions}
Let $\boldsymbol{Y}$ follow a multinomial splitting distribution $\mathcal{M}_{\Delta_{N}}( \boldsymbol{\pi}) \underset{N}{\wedge} \mathcal{L}(\psi)$. 
The criteria \ref{criter1} and \ref{criter3} hold, as a consequence from the positivity and the symmetry.  
The \gls{pmf} is given by
\begin{equation}
\label{mass:multinomial:splitting}
	P(\boldsymbol{Y}=\boldsymbol{y}) = P(\vert\boldsymbol{Y}\vert=\vert\boldsymbol{y}\vert) \binom{\vert\boldsymbol{y}\vert}{\boldsymbol{y}} \prod_{j=1}^J \pi_j^{y_j}.
\end{equation}
The expectation and covariance
are given by
\begin{align}
\operatorname{E}\left(\boldsymbol{Y}\right) & = \mu_1 \cdot \boldsymbol{\pi}, \label{eq:cmultinomial:E} \\
\operatorname{Cov}\left(\boldsymbol{Y}\right) & =  \mu_1 \cdot  \textrm{diag}(\boldsymbol{\pi}) + (\mu_2-\mu_1^2) \cdot \boldsymbol{\pi} \boldsymbol{\pi}^t,  \label{eq:cmultinomial:C}
\end{align}
where $\boldsymbol{\pi}^t$ denotes the transposition of the vector $\boldsymbol{\pi}$ and $\mu_k$ denotes the factorial moments of the sum distribution for $k=1, 2$. 
Moreover, using the \gls{pgf} formula of the singular multinomial distribution \citep{JKB97} 
we obtain the \gls{pgf} of multinomial splitting distributions as
\begin{equation}
\label{pgf:multinomial:splitting}
    G_{\boldsymbol{Y} } (\boldsymbol{s}) = \operatorname{E}\left\lbrace (\boldsymbol{\pi}^t \boldsymbol{s} )^{\vert \boldsymbol{Y} \vert} \right\rbrace = G_{\vert \boldsymbol{Y} \vert} \left( \boldsymbol{\pi}^t \boldsymbol{s} \right),
\end{equation}
where $G_{\vert \boldsymbol{Y} \vert} $ denote the \gls{pgf} of the sum distribution.
The graphical model is characterized by the following property. 
\begin{property}\citep{PF17}\label{graph_multinomial}
The minimal graphical model for a multinomial splitting distribution is empty if the sum follows a Poisson distribution and is complete otherwise.
\end{property}
Therefore, all multinomial splitting distribution are \textit{sensu stricto} multivariate distributions (the criteria \ref{criter2} holds) except when the sum follows a Poisson distribution.
As a consequence from additivity, Theorem~\ref{marginal_csd} holds and yields the marginals distributions in Corollary \label{marginal_multinomial}.
\begin{corollary}
    \label{marginal_multinomial}
	Let $\boldsymbol{Y}$ follow a multinomial splitting distribution, $\boldsymbol{Y} \sim \mathcal{M}_{\Delta_{N}}( \boldsymbol{\pi}) \underset{N}{\wedge} \mathcal{L}(\boldsymbol{\psi})$.
	Then, the marginals follow the binomial damage distribution $Y_j \sim \mathcal{B}_N(\pi_j) \underset{N}{\wedge} \mathcal{L}(\boldsymbol{\psi})$.
    Moreover, for $y\in\mathbb{N}$
    \begin{equation}
        \label{binomial-operator}
        P\left(Y_j = y\right) = \frac{\pi_j^y}{y!} G_{\vert \boldsymbol{Y} \vert}^{(y)}(1-\pi_j), 
    \end{equation} 
    where $G_{\vert \boldsymbol{Y} \vert}^{(y)}$ denotes the $y$-th derivative of the \gls{pgf} of the sum distribution $\mathcal{L}(\boldsymbol{\psi})$.
\end{corollary}
Using equation \eqref{binomial-operator}, it is easy to study the action of the binomial distribution among the set of univariate distribution.
Assume that $\mathcal{L}(\psi)$ is a power series distribution denoted by $PSD\{g(\alpha)\}$. 
It can be seen by identifiability that the resulting distribution of $\boldsymbol{Y}$ is exactly the multivariate sum-symmetric power series distribution (MSSPSD) introduced by \citet{Pat68} with
$
 	\mathcal{M}_{\Delta_N}(\boldsymbol{\pi}) \underset{N}{\wedge} PSD\{g(\alpha) \} = \textnormal{MSSPSD}\{ \alpha \cdot \boldsymbol{\pi} \}
$.
The non-singular multinomial distribution, the negative multinomial distribution and the multivariate logarithmic series distribution are thereby encompassed in multinomial splitting distributions (see Table \ref{table:multinomial:psd}). 

Assume now that $\mathcal{L}(\psi)$ is a standard beta compound distribution. 
We obtaine three new multivariate distributions which are multivariate extensions of the non-standard beta binomial, non-standard beta negative binomial and beta Poisson distributions (see Table \ref{table:multinomial:beta_compound} for details about these three multivariate distributions and Appendix \ref{RDUD} for definitions of the non-standard beta binomial and the non-standard beta negative binomial distributions).
All the characteristics of these six multinomial splitting distributions (\gls{pmf}, expectation, covariance, \gls{pgf} and marginal distributions) have been calculated using equations \eqref{mass:multinomial:splitting}, \eqref{eq:cmultinomial:E}, \eqref{eq:cmultinomial:C}, \eqref{pgf:multinomial:splitting}, \eqref{binomial-operator} according to the sum distribution $\mathcal{L}(\psi)$.
The singular multinomial distribution belongs to the exponential family.
For any $j \in \mathcal \{1,\ldots,J\}$  the \glspl{mle} $\widehat{\pi_{j}}$ are given by the following closed-form expression:
\begin{equation}
    \label{eq:multinomial:distribution:mle}
    \widehat{\pi}_j = \frac{\sum_{\boldsymbol{y} \in \boldsymbol{\mathcal{N}}} {y_j}}
{\sum_{\boldsymbol{y} \in \boldsymbol{\mathcal{N}}} \vert \boldsymbol{y} \vert }.
\end{equation}

\subparagraph{Canonical case of binomial sum distribution.}
The case $\mathcal{L}(\psi)=\mathcal{B}_{n}(p)$ is considered as the canonical case since the binomial distribution is the univariate version of the non-singular multinomial distribution.
Usual characteristics of the multinomial splitting binomial distribution are obtained using equations \eqref{mass:multinomial:splitting}, \eqref{eq:cmultinomial:E}, \eqref{eq:cmultinomial:C}, \eqref{pgf:multinomial:splitting}, \eqref{binomial-operator} with $\mathcal{L}(\psi)=\mathcal{B}_{n}(p)$.
It should be remarked that the constraint between parameters of the singular distribution and the sum distribution described in Corollary~\ref{coro:convolution:nonsingular} disappears in this case, as for equation \eqref{msd}.
More generally the constraint disappears when the convolution is proportional.
Relation~\eqref{msd} and equation~\eqref{prop:inference:freq} allow parameters of non-singular multinomial distributions to be inferred using a parameter inference procedure for singular multinomial and binomial distributions. 
In particular, this enables us to propose joint or conditional maximum likelihood estimation of parameters $n$ and $\boldsymbol{\pi}$ of the non-singular multinomial distribution directly by combining \gls{mle} inference for the binomial distribution \citep{BD81} with the closed-form expression given in equation~\eqref{eq:multinomial:distribution:mle}.
Moreover, equation~\eqref{prop:inference:freq} allows us to generalize studies related to  parameter inference in binomial distributions \citep{Der83} to parameter inference in the non-singular multinomial distributions. If the sum dataset is 
\begin{itemize}
	\item overdispersed, the likelihood is an increasing function of $n$ and there is no finite \gls{mle} of $n$,
	\item underdispersed, the likelihood is either a decreasing function of $n$, or has an unique maximum and there is a finite \gls{mle} of $n$.
\end{itemize}
Note that some previous works considered the joint estimation of parameters $n$ and $\boldsymbol{\pi}$ for non-singular multinomial distributions \citep{San72,RDL07} but only \glspl{mle} for specific constraints were derived.

\begin{landscape}
\begin{table}
\begin{tabular}{|l||l|l|l|}
\hline
Distribution & \multicolumn{3}{c|}{$\boldsymbol{Y} \sim \mathcal{M}_{\Delta_N}(\boldsymbol{\pi}) \underset{N}{\wedge} \mathcal{L}(\psi)$} \\
\hline
$\mathcal{L}(\psi)$
& $\mathcal{B}_{n}(p)$ 
& $\mathcal{NB}(r, p)$ 
& $\mathcal{L}(p)$ 
\\
\hline
Re-parametrization
& $\mathcal{M}_{\blacktriangle_{n}}(p\cdot\boldsymbol{\pi})$
& $\mathcal{NM}(r, p\cdot\boldsymbol{\pi})$
& $\mathcal{ML}(p\cdot\boldsymbol{\pi})$
\\
\hline
Supp$(\boldsymbol{Y})$ 
& $\blacktriangle_n$ 
& $\mathbb{N}^J$ 
& $\mathbb{N}^J \setminus (0,\ldots,0)$ 
\\
\hline
$P(\boldsymbol{Y}=\boldsymbol{y})$ 
& $\binom{n}{\boldsymbol{y}}  \left(1 - p \right)^{n - \vert\boldsymbol{y}\vert } \prod_{j=1}^J (p\pi_j)^{y_j}$  
& $\binom{\vert\boldsymbol{y}\vert+r-1}{\boldsymbol{y}}  \left(1 - p \right)^{r} \prod_{j=1}^J (p\pi_j)^{y_j}$
& $\binom{\vert\boldsymbol{y}\vert}{\boldsymbol{y}} \frac{-1}{\vert\boldsymbol{y}\vert\ln(1-p)} \prod_{j=1}^J \pi_j^{y_j}$
\\
\hline
$\operatorname{E}(\boldsymbol{Y})$
& $n p \cdot \boldsymbol{\pi}$
& $r \frac{p}{1-p} \cdot \boldsymbol{\pi}$
& $\frac{-p}{(1-p)\ln (1-p)} \cdot \boldsymbol{\pi}$
\\
\hline
$\operatorname{Cov}(\boldsymbol{Y})$
& $n p \cdot \left\lbrace \textrm{diag}(\boldsymbol{\pi}) - p \cdot \boldsymbol{\pi}\boldsymbol{\pi}^t  \right\rbrace$
& $r \frac{p}{1-p} \cdot \left\lbrace \textrm{diag}(\boldsymbol{\pi}) + \frac{p}{1-p} \cdot \boldsymbol{\pi}\boldsymbol{\pi}^t  \right\rbrace$
& $\frac{-p}{(1-p)\ln (1-p)} \cdot \left\lbrace \textrm{diag}(\boldsymbol{\pi}) + \frac{p\{1-\ln (1-p)\}}{(1-p)\ln (1-p)} \cdot \boldsymbol{\pi}\boldsymbol{\pi}^t  \right\rbrace$
\\
\hline
$G_{\boldsymbol{Y}}(\boldsymbol{s})$
& $\left( 1-p+p\; \boldsymbol{\pi}^t \boldsymbol{s} \right)^{n}$
& $\left( \frac{1-p}{1-p \boldsymbol{\pi}^t \boldsymbol{s} }\right)^{r}$
& $\frac{\ln (1-p \cdot \boldsymbol{\pi}^t \boldsymbol{s})}{\ln(1-p)}$
\\
\hline
Marginals
& $Y_j\sim\mathcal{B}_{n}(p\pi_j)$ \footnote{This is a natural multivariate extension of $\mathcal{L}(\boldsymbol{\psi})$}\saveFN\sfn
& $Y_j\sim\mathcal{NB}(r, p\pi_j)$\useFN\sfn\ 
& $Y_j\sim\mathcal{L}\left(p^{\prime},\omega\right)$ \footnote{where $\omega=\frac{\ln (1-p+p\pi_j)}{\ln (1-p)}$ and $p^{\prime}=\frac{p\pi_j}{1-p+p\pi_j}$}
\\
\hline
\end{tabular}
\caption{\label{table:multinomial:psd}Usual characteristics 
of multinomial splitting binomial, negative binomial and logarithmic series distribution.}
\end{table}

\begin{table}
\begin{tabular}{|l||l|l|l|}
\hline
Distribution & \multicolumn{3}{c|}{$\boldsymbol{Y} \sim \mathcal{M}_{\Delta_N}(\boldsymbol{\pi}) \underset{N}{\wedge} \mathcal{L}(\psi)$} \\
\hline
$\mathcal{L}(\psi)$
& $\beta\mathcal{B}_{n}(a,b)$ 
& $\beta\mathcal{NB}(r,a,b)$ 
& $\beta_{\lambda}\mathcal{P}(a,b)$ 
\\
\hline
Supp$(\boldsymbol{Y})$ 
& $\blacktriangle_n$ 
& $\mathbb{N}^J$ 
& $\mathbb{N}^J$ 
\\
\hline
$P(\boldsymbol{Y}=\boldsymbol{y})$ 
& $\binom{n}{\boldsymbol{y}}  \frac{B(\vert\boldsymbol{y}\vert+a,n-\vert\boldsymbol{y}\vert+b)}{B(\alpha,b)}  \prod_{j=1}^J \pi_j^{y_j} $
& $\binom{\vert\boldsymbol{y}\vert+r-1}{\boldsymbol{y},r-1}  \frac{B(r+a,\vert\boldsymbol{y}\vert+b)}{B(a,b)}  \prod_{j=1}^J \pi_j^{y_j}$
& $\frac{(a)_{\vert\boldsymbol{y}\vert} \lambda^{\vert\boldsymbol{y}\vert}}{(a+b)_{\vert\boldsymbol{y}\vert}} {}_1F_1(a+\vert\boldsymbol{y}\vert,a+b+\vert\boldsymbol{y}\vert;-\lambda) \prod_{j=1}^J \frac{\pi_j^{y_j}}{y_j!}$
\\
\hline
$\operatorname{E}(\boldsymbol{Y})$
& $n \frac{a}{a+b} \cdot \boldsymbol{\pi}$
& $r \frac{b}{a-1} \cdot \boldsymbol{\pi}$ \footnote{if $\alpha>1$ and not defined otherwise}
& $\lambda \frac{a}{a+b} \cdot \boldsymbol{\pi}$
\\
\hline
$\operatorname{Cov}(\boldsymbol{Y})$
& $n \frac{a}{a+b} \cdot \left\lbrace \textrm{diag}(\boldsymbol{\pi}) + \frac{b(n-1)-a(a+b+1)}{(a+b)(a+b+1)} \cdot \boldsymbol{\pi}\boldsymbol{\pi}^t  \right\rbrace$
& $r \frac{b}{a-1} \cdot \left\lbrace \textrm{diag}(\boldsymbol{\pi}) + \frac{a(b+r+1)+r(b-1)-b-1}{(a-1)(a-2)} \cdot \boldsymbol{\pi}\boldsymbol{\pi}^t  \right\rbrace$ \footnote{if $\alpha>2$ and not defined otherwise}
& $\lambda \frac{a}{a+b} \cdot \left\lbrace \textrm{diag}(\boldsymbol{\pi}) + \lambda \frac{b}{(a+b)(a+b+1)} \cdot \boldsymbol{\pi}\boldsymbol{\pi}^t  \right\rbrace$
\\
\hline
$G_{\boldsymbol{Y}}(\boldsymbol{s})$
& ${}_2F_1\{(-n, a); \alpha+b; 1-\boldsymbol{\pi}^t \boldsymbol{s}\}$
& $\frac{a_{(b)}}{(r+a)_{(b)}} {}_2F_1\left\lbrace (r, b);r+a+b; \boldsymbol{\pi}^t \boldsymbol{s} \right\rbrace$
& ${}_1F_1\left\lbrace a; a+b; \lambda(\boldsymbol{\pi}^t \boldsymbol{s}-1) \right\rbrace$
\\
\hline
Marginals
& $Y_j\sim\beta_{\pi_j}\mathcal{B}_{n}(a,b)$
& $Y_j\sim\beta_{\pi_j}\mathcal{NB}(r,a,b)$
& $Y_j\sim\beta_{\pi_j\lambda}\mathcal{P}(a,b)$
\\
\hline
\end{tabular}
\caption{\label{table:multinomial:beta_compound}Usual characteristics of multinomial splitting beta binomial, beta negative binomial and beta Poisson distribution.}
\end{table}
\end{landscape}

\section{Dirichlet Multinomial Splitting Models}
\label{section:dirichlet:multinomial}
In this section the Dirichlet multinomial distribution is introduced as a positive and additive convolution distribution.
Then, the general case of Dirichlet multinomial splitting distributions is studied.
For six specific sum distributions, the usual characteristics of Dirichlet multinomial splitting distributions are described in Tables \ref{table:dirichlet:multinomial:beta_binomial}, \ref{table:dirichlet:multinomial:beta_negative_binomial}, 
\ref{table:dirichlet:multinomial:beta_poisson},
and \ref{table:dirichlet:multinomial:psd}.
Finally, the canonical case of beta binomial sum distribution is detailed, with particular emphasis on parameter inference.

\subparagraph{Dirichlet multinomial distribution}
Let $a_{\theta}(y)=\binom{y+\theta-1}{y}$ be the parametric sequence that characterizes the Dirichlet multinomial distribution as a convolution distribution.
It is positive since $\binom{y+\theta-1}{y}>0$ for all $\theta\in\Theta=(0,\infty)$ and all $y\in\mathbb{N}$. 
It is additive, as a consequence from the convolution identity of Hagen and Rothe: $\binom{n+\theta+\gamma-1}{n}=\sum_{y=0}^n \binom{y+\theta-1}{y} \binom{n-y+\gamma-1}{n-y}$.
It implies, by induction on $n$, that the normalizing constant is  $c_{\boldsymbol{\theta}}(n)=a_{\vert\boldsymbol{\theta}\vert}(n)=\binom{n+\vert\boldsymbol{\theta}\vert-1}{n}$. 
In order to respect the usual notation, parameter $\boldsymbol{\alpha}$ will be used instead of $\boldsymbol{\theta}$, and thus the Dirichlet multinomial distribution will be denoted by $\mathcal{DM}_{\Delta_n}(\boldsymbol{\alpha})$ with $n\in\mathbb{N}$ and $\boldsymbol{\alpha}\in(0,\infty)^J$.
The non-singular Dirichlet multinomial distribution will be denoted by $\mathcal{DM}_{\blacktriangle_{n}}( \boldsymbol{\alpha}, b)$ with $b\in(0,\infty)$. 
The beta binomial distribution will be denoted by $\beta \mathcal{B}_n(a,b)$ with $(a,b)\in(0,\infty)^2$. 
For similar reasons as in the multinomial case, the non-singular Dirichlet multinomial distribution should be considered as the natural extension of the beta binomial distribution, rather than the singular one.
Note that a Dirichlet multinomial distribution is exactly a multivariate negative hypergeometric distribution. 

\subparagraph{Dirichlet multinomial splitting distributions}
In this paragraph let $\boldsymbol{Y}$ follow a Dirichlet multinomial splitting distribution $\mathcal{DM}_{\Delta_{N}}( \boldsymbol{\alpha}) \underset{N}{\wedge} \mathcal{L}(\psi)$. 
Criteria \ref{criter1} and \ref{criter3} hold, as a consequence of positivity and symmetry.  
The \gls{pmf} is given by
\begin{equation}
\label{mass:dirichlet:multinomial:splitting}
	P(\boldsymbol{Y}=\boldsymbol{y}) =\frac{ P(\vert\boldsymbol{Y}\vert=\vert\boldsymbol{y}\vert)  }{ \binom{ \vert\boldsymbol{y}\vert +\vert \boldsymbol{\alpha} \vert -1}{\vert\boldsymbol{y}\vert} }   \prod_{j=1}^J \binom{y_j+\alpha_j-1}{y_j}.
\end{equation}
The expectation and covariance are given by
\begin{align}
	\operatorname{E}\left(\boldsymbol{Y}\right) 
    	& = \frac{\mu_1}{\vert\boldsymbol{\alpha}\vert} \cdot \boldsymbol{\alpha}, \label{eq:dirichlet:cmultinomial:E} \\
	\operatorname{Cov}\left(\boldsymbol{Y}\right) 
    	& =  \frac{1}{\vert\boldsymbol{\alpha}\vert(\vert\boldsymbol{\alpha}\vert+1)} \cdot \left[ \left\lbrace (\vert\boldsymbol{\alpha}\vert+1)\mu_1 +\mu_2 \right\rbrace \cdot \textnormal{diag}(\boldsymbol{\alpha}) +  \left\lbrace \mu_2 - \frac{\vert\boldsymbol{\alpha}\vert+1}{\vert\boldsymbol{\alpha}\vert}\mu_1^2 \right\rbrace \cdot \boldsymbol{\alpha}\boldsymbol{\alpha}^t \right]. \label{eq:dirichlet:cmultinomial:C}
\end{align}
The \gls{pgf} is given by
\begin{equation}\label{pgf:dirichlet:multinomial:splitting}
	G_{\boldsymbol{Y}}(\boldsymbol{s}) = \sum_{\boldsymbol{y} \in \mathbb{N}^J} \Gamma\left( \vert \boldsymbol{y} \vert +1 \right)  P\left( \vert \boldsymbol{Y} \vert = \vert \boldsymbol{y} \vert \right) \frac{\prod_{j=1}^J(\alpha_j)_{y_j}}{(\vert\boldsymbol{\alpha}\vert)_{\vert\boldsymbol{y}\vert}} \prod_{j=1}^J \frac{s_j^{y_j}}{y_j!}.
\end{equation}
The graphical model is characterized by the following property. 
\begin{property}\citep{PF17}\label{graph:dirichlet:multinomial}
Let $\boldsymbol{Y}$ follow a Dirichlet multinomial splitting distribution with parameter $\boldsymbol{\alpha} \in (0,\infty)^J$.
The minimal graphical model for $\boldsymbol{Y}$ is empty if the sum follows a negative binomial distribution with parameters $(\vert \boldsymbol{\alpha} \vert, p)$ for some $p \in (0,1)$ and is complete otherwise.
\end{property}
Therefore, all Dirichlet multinomial splitting distribution are \textit{senso stricto} multivariate distributions except when the sum follows a negative binomial distribution $\mathcal{NB}(r,p)$ with the specific constraint $r = \vert \boldsymbol{\alpha} \vert$. 
Finally, due to additivity, Property~\ref{marginal_csd} can be applied to describe the marginal distributions.
\begin{corollary}
    \label{marginal:dirichlet:multinomial}
	Let $\boldsymbol{Y}$ follow a Dirichlet multinomial splitting distribution, $\boldsymbol{Y} \sim \mathcal{DM}_{\Delta_{N}}( \boldsymbol{\alpha}) \underset{N}{\wedge} \mathcal{L}(\boldsymbol{\psi})$ with $\boldsymbol{\alpha} \in (0,\infty)^J$.
    Then, the marginals follow the beta-binomial damage distribution $Y_j \sim \beta\mathcal{B}_N(\alpha_j, \vert \boldsymbol{\alpha}_{-j} \vert ) \underset{N}{\wedge} \mathcal{L}(\boldsymbol{\psi})$.
\end{corollary}
Using the Fubini theorem, it can be shown that 
\begin{equation} \label{beta-binomial-operator}
Y_j \sim \left\lbrace \mathcal{B}_N(\pi) \underset{N}{\wedge} \mathcal{L}(\boldsymbol{\psi}) \right\rbrace \underset{\pi}{\wedge} \beta(\alpha_j, \vert \boldsymbol{\alpha}_{-j} \vert ),
\end{equation}
since $N$ and $\pi$ are independent latent variables.
Therefore, results previously obtained for the binomial damage distributions can be used to describe the beta-binomial damage distributions.
Assume that $\mathcal{L}(\psi)$ is a standard beta compound distribution.
Four new and two already known multivariate distributions are obtained  or recovered (see Tables~\ref{table:dirichlet:multinomial:beta_binomial},  \ref{table:dirichlet:multinomial:beta_negative_binomial} and \ref{table:dirichlet:multinomial:beta_poisson}). 
In particular, natural multivariate extensions of three beta compound distributions are described. 
The non-singular Dirichlet multinomial is recovered when $\mathcal{L}(\psi)=\beta\mathcal{B}_n(a,b)$ with the specific constraint $a=\vert\boldsymbol{\alpha}\vert$.
The multivariate generalized waring distribution, introduced by \citet{Xek86}, is recovered when $\mathcal{L}(\psi)=\beta\mathcal{NB}(r,a,b)$ with the specific constraint $r=\vert\boldsymbol{\alpha}\vert$.
Finally, a multivariate extension of the beta Poisson distribution is proposed when $\mathcal{L}(\psi)=\beta_{\lambda}\mathcal{P}(a,b)$ with the specific constraint $a=\vert\boldsymbol{\alpha}\vert$. 

Assume now that $\mathcal{L}(\psi)$ is a power series distributions leading to three new multivariate extensions (see Table \ref{table:dirichlet:multinomial:psd}).
Remark that several multivariate extensions of the same univariate distribution could be defined.
For instance the multinomial splitting beta binomial distribution $\mathcal{M}_{\Delta_N}(\boldsymbol{\pi})\underset{N}{\wedge}\beta\mathcal{B}_n(a,b)$ and the Dirichlet multinomial splitting binomial distribution $\mathcal{DM}_{\Delta_N}(\boldsymbol{\alpha})\underset{N}{\wedge}\mathcal{B}_n(p)$ are two multivariate extensions of the non-standard beta binomial distribution (see Tables \ref{table:multinomial:beta_compound} and \ref{table:dirichlet:multinomial:psd}).
Note that the singular Dirichlet multinomial distribution does not belongs to the exponential family.
Either if $\left\vert \boldsymbol{\alpha} \right\vert$ is known or not, \glspl{mle} $\widehat{\alpha_{j}}$ can be computed using various iterative methods \citep{Min00, Skl14}.

\subparagraph{Canonical case of beta binomial sum distribution}
The case $\mathcal{L}(\psi)=\beta\mathcal{B}_{n}(a,b)$ is considered as the canonical case since the beta binomial distribution is the univariate version of the non-singular Dirichlet multinomial distribution.
Usual characteristics of the Dirichlet multinomial splitting beta binomial distribution 
are derived from equations \eqref{mass:dirichlet:multinomial:splitting}, \eqref{eq:dirichlet:cmultinomial:E}, \eqref{eq:dirichlet:cmultinomial:C}, \eqref{pgf:dirichlet:multinomial:splitting} and \eqref{beta-binomial-operator} with $\mathcal{L}(\psi)=\beta\mathcal{B}_{n}(a,b)$.
According to Corollary~\ref{coro:convolution:nonsingular}, the Dirichlet multinomial splitting beta binomial distribution with the specific constraint $a=\vert\boldsymbol{\alpha}\vert$ is exactly the non-singular Dirichlet multinomial distribution: 
\begin{equation}\label{dmsd}  
	\mathcal{DM}_{\Delta_N}\left( \boldsymbol{\alpha}\right) \underset{N}{\wedge} \beta\mathcal{B}_{n}\left( \vert \boldsymbol{\alpha} \vert , b  \right) = \mathcal{DM}_{\blacktriangle_{n}}\left( \boldsymbol{\alpha}, b \right).
\end{equation}
The constraint $a=\vert\boldsymbol{\alpha}\vert$ has to be taken into account in the inference procedure, either on the singular distribution or on the sum distribution.
We propose to use the first alternative since the inference procedure of a constrained Dirichlet multinomial distribution (i.e., with a fixed sum $\vert\boldsymbol{\alpha}\vert$) has already been proposed by \cite{Min00}.
The sum distribution $\beta\mathcal{B}_{n}\left(a, b\right)$ can then be estimated without constraint on parameters $a$ or $b$ (see Table~\ref{table:dud}).
Note that, if no constraint between parameters of singular and sum distributions is assumed then the inference procedure is straightforward, since it can be separated into two independent procedures.
The resulting splitting distribution is more general, including the non-singular Dirichlet multinomial distribution as a special case.
As a consequence from equation \eqref{beta-binomial-operator}, the marginals follow beta square binomial distributions $\beta^2\mathcal{B}_n(\alpha_j,\vert\boldsymbol{\alpha}_{-j}\vert,a,b)$ and $\beta\mathcal{B}_n(\alpha_j,\vert\boldsymbol{\alpha}_{-j}\vert+b)$ when the constraint $a=\vert\boldsymbol{\alpha}\vert$ is assumed (see Appendices \ref{RCUD} and \ref{RDUD} for definition of beta square distribution and beta square compound distributions).
Identifying these two distributions we obtain a property about the product of two independent beta distributions.
\begin{property}
For $(a,b,c)\in(0,\infty)^3$, let $X\sim \beta (a,b)$ and $Y\sim \beta(a+b,c)$ be two independent random variables. Then $XY\sim \beta(a,b+c)$.
\end{property}
This result can be extended by induction for a product of $n$ independent beta distributions. 

\begin{landscape}
\begin{table}
\begin{tabular}{|l||l|l|}
\hline
Distribution & \multicolumn{2}{c|}{$\boldsymbol{Y}\sim\mathcal{DM}_{\Delta_N}(\boldsymbol{\alpha})\underset{N}{\wedge}\beta\mathcal{B}_{n}(a,b)$ } \\
\hline
Constraint
& no constraint
& $a=\vert\boldsymbol{\alpha}\vert$ 
\\
\hline
Re-parametrization
& 
& $\mathcal{DM}_{\blacktriangle_{n}}(\boldsymbol{\alpha},b)$
\\
\hline
Supp$(\boldsymbol{Y})$ 
& $\blacktriangle_n$ 
& $\blacktriangle_n$
\\
\hline
$P(\boldsymbol{Y}=\boldsymbol{y})$ 
& $\binom{n}{\vert\boldsymbol{y}\vert} \frac{B(a+\vert\boldsymbol{y}\vert,b+n-\vert\boldsymbol{y}\vert)}{B(a,b)} \frac{  \prod_{j=1}^J \binom{y_j+\alpha_j-1}{y_j} }{ \binom{ n +\vert \boldsymbol{\alpha} \vert -1}{n}}$
& $\binom{n - \vert \boldsymbol{y} \vert +b-1}{n - \vert \boldsymbol{y} \vert } \frac{  \prod_{j=1}^J \binom{y_j+\alpha_j-1}{y_j} }{ \binom{ n +\vert \boldsymbol{\alpha} \vert + b -1}{n} }$
\\
\hline
$\operatorname{E}(\boldsymbol{Y})$
& $\frac{n a}{\vert\boldsymbol{\alpha}\vert(a+b)} \cdot \boldsymbol{\alpha}$
& $\frac{n\vert\boldsymbol{\alpha}\vert}{\vert\boldsymbol{\alpha}\vert(\vert\boldsymbol{\alpha}\vert+b)} \cdot \boldsymbol{\alpha}$
\\
\hline
$\operatorname{Cov}(\boldsymbol{Y})$
& 
\begin{tabular}[t]{l} 
$\frac{na}{\vert\boldsymbol{\alpha}\vert(\vert\boldsymbol{\alpha}\vert+1)(a+b)} \cdot\bigg[ \left\lbrace \frac{b(a+b+n)}{(a+b)(a+b+1)} + \frac{na}{a+b}  + \vert\boldsymbol{\alpha}\vert \right\rbrace \cdot\textrm{diag}(\boldsymbol{\alpha})$  \\ 
$ + \left\lbrace \frac{b(a+b+n)}{(a+b)(a+b+1)} - \frac{na}{\vert\boldsymbol{\alpha}\vert(a+b)} -1 \right\rbrace \cdot\boldsymbol{\alpha}\boldsymbol{\alpha}^t \bigg]$ 
\end{tabular}
&
\begin{tabular}[t]{l} 
$\frac{n}{(\vert\boldsymbol{\alpha}\vert+1)(\vert\boldsymbol{\alpha}\vert+b)} \cdot\bigg[ \left\lbrace  \frac{b(\vert\boldsymbol{\alpha}\vert+b+n)}{(\vert\boldsymbol{\alpha}\vert+b)(\vert\boldsymbol{\alpha}\vert+b+1)} + \frac{n\vert\boldsymbol{\alpha}\vert}{\vert\boldsymbol{\alpha}\vert+b}  + \vert\boldsymbol{\alpha}\vert \right\rbrace \cdot\textrm{diag}(\boldsymbol{\alpha})$  \\ 
$ + \left\lbrace \frac{b(\vert\boldsymbol{\alpha}\vert+b+n)}{(\vert\boldsymbol{\alpha}\vert+b)(\vert\boldsymbol{\alpha}\vert+b+1)} - \frac{n}{\vert\boldsymbol{\alpha}\vert+b} -1 \right\rbrace \cdot\boldsymbol{\alpha}\boldsymbol{\alpha}^t \bigg]$ 
\end{tabular}
\\
\hline
$G_{\boldsymbol{Y}}(\boldsymbol{s})$
& 
$\frac{(b)_{n}}{(a+b)_{n}} {}_{2}^JF_2\{(-n,a);\boldsymbol{\alpha};(-b-n+1,\vert\boldsymbol{\alpha}\vert);\boldsymbol{s}\}$
& $\frac{(b)_{n}}{(\vert\boldsymbol{\alpha}\vert+b)_{n}} {}_{1}^JF_1(-n;\boldsymbol{\alpha};-b-n+1;\boldsymbol{s})$
\\
\hline
Marginals
& $Y_j\sim\beta^2\mathcal{B}_{n}\left(\alpha_j, \vert \boldsymbol{\alpha}_{-j}\vert, a, b \right)$
& $Y_j\sim\beta\mathcal{B}_{n}\left(\alpha_j, \vert\boldsymbol{\alpha}_{-j}\vert + b\right)$\useFN\sfn\
\\
\hline
\end{tabular}
\caption{\label{table:dirichlet:multinomial:beta_binomial}Usual characteristics 
of Dirichlet multinomial splitting standard beta binomial distribution respectively without constraint and with $a=\vert\boldsymbol{\alpha}\vert$ .}
\end{table}

\begin{table}
\begin{tabular}{|l||l|l|}
\hline
Distribution  & \multicolumn{2}{c|}{$\boldsymbol{Y}\sim\mathcal{DM}_{\Delta_N}(\boldsymbol{\alpha})\underset{N}{\wedge}\beta\mathcal{NB}(r,a,b)$} \\
\hline
Constraint
& no constraint
& $r=\vert\boldsymbol{\alpha}\vert$ 
\\
\hline
Re-parametrization
& 
& MGWD$(b,\boldsymbol{\alpha},a)$
\\
\hline
Supp$(\boldsymbol{Y})$ 
& $\mathbb{N}^J$
& $\mathbb{N}^J$ 
\\
\hline
$P(\boldsymbol{Y}=\boldsymbol{y})$ 
& $\frac{(a)_{r}}{(a+b)_{r}} \frac{(r)_{\vert\boldsymbol{y}\vert}(b)_{\vert\boldsymbol{y}\vert}}{(r+a+b)_{\vert\boldsymbol{y}\vert}(\vert\boldsymbol{\alpha}\vert)_{\vert\boldsymbol{y}\vert}}\prod_{j=1}^J\frac{(\alpha_j)_{y_j}}{y_j!}$
& $\frac{(a)_{\vert\boldsymbol{\alpha}\vert}}{(a+b)_{\vert\boldsymbol{\alpha}\vert}} \frac{(b)_{\vert\boldsymbol{y}\vert}}{(\vert\boldsymbol{\alpha}\vert+a+b)_{\vert\boldsymbol{y}\vert}}\prod_{j=1}^J\frac{(\alpha_j)_{y_j}}{y_j!}$
\\
\hline
$\operatorname{E}(\boldsymbol{Y})$
& $\frac{rb}{\vert\boldsymbol{\alpha}\vert(\vert\boldsymbol{\alpha}\vert+1)(a-1)} \boldsymbol{\alpha} $
& $\frac{b}{(\vert\boldsymbol{\alpha}\vert+1)(a-1)} \boldsymbol{\alpha} $
\\
\hline
$\operatorname{Cov}(\boldsymbol{Y})$
& 
\begin{tabular}[t]{l} 
$\frac{rb}{\vert\boldsymbol{\alpha}\vert(\vert\boldsymbol{\alpha}\vert+1)(a-1)} 
\cdot\bigg[ \left\lbrace  \frac{(r+a-1)(a+b-1)}{(a-1)(a-2)} + \frac{rb}{a-1} + \vert\boldsymbol{\alpha}\vert \right\rbrace \cdot\textrm{diag}(\boldsymbol{\alpha})$  \\ 
$ + \left\lbrace \frac{(r+a-1)(a+b-1)}{(a-1)(a-2)} - \frac{rb}{\vert\boldsymbol{\alpha}\vert(a-1)} -1 \right\rbrace \cdot\boldsymbol{\alpha}\boldsymbol{\alpha}^t \bigg]$ \footnote{defined if $a>1$}
\end{tabular}
&
\begin{tabular}[t]{l} 
$\frac{b}{(\vert\boldsymbol{\alpha}\vert+1)(a-1)} 
\cdot\bigg[ \left\lbrace \frac{(\vert\boldsymbol{\alpha}\vert+a-1)(a+b-1)}{(a-1)(a-2)} + \frac{\vert\boldsymbol{\alpha}\vert b}{a-1} + \vert\boldsymbol{\alpha}\vert \right\rbrace \cdot\textrm{diag}(\boldsymbol{\alpha})$  \\ 
$ + \left\lbrace \frac{(\vert\boldsymbol{\alpha}\vert+a-1)(a+b-1)}{(a-1)(a-2)} - \frac{b}{a-1} -1 \right\rbrace \cdot\boldsymbol{\alpha}\boldsymbol{\alpha}^t \bigg]$ \footnote{defined if $a>2$}
\end{tabular}
\\
\hline
$G_{\boldsymbol{Y}}(\boldsymbol{s})$
& $\frac{(a)_{r}}{(a+b)_{r}} {}_{2}^JF_2\{(r,b);\boldsymbol{\alpha};(r+a+b,\vert\boldsymbol{\alpha}\vert);\boldsymbol{s}\}$
& $\frac{(a)_{\vert\boldsymbol{\alpha}\vert}}{(a+b)_{\vert\boldsymbol{\alpha}\vert}} {}_{1}^JF_1(b;\boldsymbol{\alpha};\vert\boldsymbol{\alpha}\vert+a+b;\boldsymbol{s})$
\\
\hline
Marginals
& $Y_j\sim\beta^2\mathcal{NB}\left(r,\alpha_j, \vert \boldsymbol{\alpha}_{-j}\vert, a, b \right)$
& $Y_j\sim\beta\mathcal{NB}\left(\alpha_j, a, b\right)$\useFN\sfn\
\\
\hline 
\end{tabular}
\caption{\label{table:dirichlet:multinomial:beta_negative_binomial}Usual characteristics 
of Dirichlet multinomial splitting standard beta negative binomial distribution respectively without constraint and with $r=\vert\boldsymbol{\alpha}\vert$.}
\end{table}
\end{landscape}

\begin{landscape}
\begin{table}
\begin{tabular}{|l||l|l|}
\hline
Distribution & \multicolumn{2}{c|}{$\boldsymbol{Y}\sim\mathcal{DM}_{\Delta_N}(\boldsymbol{\alpha})\underset{N}{\wedge}\beta_{\lambda}\mathcal{P}(a,b)$ } \\
\hline
Constraint
& no constraint
& $a=\vert\boldsymbol{\alpha}\vert$ 
\\
\hline
Supp$(\boldsymbol{Y})$ 
& $\mathbb{N}^J$ 
& $\mathbb{N}^J$
\\
\hline
$P(\boldsymbol{Y}=\boldsymbol{y})$ 
& $\frac{(a)_{\vert\boldsymbol{y}\vert}\lambda^{\vert\boldsymbol{y}\vert}}{(a+b)_{\vert\boldsymbol{y}\vert}(\vert\boldsymbol{\alpha}\vert)_{\vert\boldsymbol{y}\vert}}\prod_{j=1}^J\frac{(\alpha_j)_{y_j}}{y_j!} {}_1F_1(a+\vert\boldsymbol{y}\vert;a+b+\vert\boldsymbol{y}\vert;-\lambda)$
& $\frac{\lambda^{\vert\boldsymbol{y}\vert}}{(a+b)_{\vert\boldsymbol{y}\vert}}\prod_{j=1}^J\frac{(\alpha_j)_{y_j}}{y_j!} {}_1F_1(a+\vert\boldsymbol{y}\vert;a+b+\vert\boldsymbol{y}\vert;-\lambda)$
\\
\hline
$\operatorname{E}(\boldsymbol{Y})$
& $\frac{\lambda a}{\vert\boldsymbol{\alpha}\vert(a+b)} \cdot \boldsymbol{\alpha}$
& $\frac{\lambda}{\vert\boldsymbol{\alpha}\vert+b} \cdot \boldsymbol{\alpha}$
\\
\hline
$\operatorname{Cov}(\boldsymbol{Y})$
& 
\begin{tabular}[t]{l} 
$\frac{\lambda a}{\vert\boldsymbol{\alpha}\vert(\vert\boldsymbol{\alpha}\vert+1)(a+b)} \cdot\bigg[ \left\lbrace \frac{\lambda b}{(a+b)(a+b+1)} + \frac{\lambda a}{a+b}  + \vert\boldsymbol{\alpha}\vert +1 \right\rbrace \cdot\textrm{diag}(\boldsymbol{\alpha})$  \\ 
$ + \left\lbrace \frac{\lambda b}{(a+b)(a+b+1)} - \frac{\lambda a}{\vert\boldsymbol{\alpha}\vert(a+b)} \right\rbrace \cdot\boldsymbol{\alpha}\boldsymbol{\alpha}^t \bigg]$ 
\end{tabular}
&
\begin{tabular}[t]{l} 
$\frac{\lambda}{(\vert\boldsymbol{\alpha}\vert+1)(a+b)} \cdot\bigg[ \left\lbrace \frac{\lambda b}{(\vert\boldsymbol{\alpha}\vert+b)(\vert\boldsymbol{\alpha}\vert+b+1)} + \frac{\lambda \vert\boldsymbol{\alpha}\vert}{\vert\boldsymbol{\alpha}\vert+b}  + \vert\boldsymbol{\alpha}\vert +1 \right\rbrace \cdot\textrm{diag}(\boldsymbol{\alpha})$  \\ 
$ + \left\lbrace \frac{\lambda b}{(\vert\boldsymbol{\alpha}\vert+b)(\vert\boldsymbol{\alpha}\vert+b+1)} - \frac{\lambda}{a+b} \right\rbrace \cdot\boldsymbol{\alpha}\boldsymbol{\alpha}^t \bigg]$ 
\end{tabular}
\\
\hline
$G_{\boldsymbol{Y}}(\boldsymbol{s})$
& $\sum_{\boldsymbol{y}\in\mathbb{N}^J}\sum_{k\in\mathbb{N}} \frac{(a)_{\vert\boldsymbol{y}\vert+k}\prod_{j=1}^J(\alpha_j)_{y_j}}{(a+b)_{\vert\boldsymbol{y}\vert+k}(\vert\boldsymbol{\alpha}\vert)_{\vert\boldsymbol{y}\vert}}\frac{(-\lambda)^k}{k!}\prod_{j\in\mathcal{J}}\frac{(\lambda s_j)^{y_j}}{y_j!}$
& ${}^{J+1}_{\;\;\;\;0}F_1\{(\boldsymbol{\alpha},a+\vert\boldsymbol{y}\vert);a+b;(\lambda\cdot\boldsymbol{s},-\lambda)\}$
\\
\hline
Marginals
& $Y_j\sim\beta^2_{\lambda}\mathcal{P}\left(\alpha_j, \vert \boldsymbol{\alpha}_{-j}\vert, a, b \right)$
& $Y_j\sim\beta_{\lambda}\mathcal{P}\left(\alpha_j, \vert\boldsymbol{\alpha}_{-j}\vert + b\right)$\useFN\sfn\
\\
\hline
\end{tabular}
\caption{\label{table:dirichlet:multinomial:beta_poisson}Usual characteristics 
of Dirichlet multinomial splitting beta Poisson distribution respectively without constraint and with $a=\vert\boldsymbol{\alpha}\vert$.}
\end{table}

\begin{table}
\begin{tabular}{|l||l|l|l|}
\hline
Distribution & \multicolumn{3}{c|}{$\boldsymbol{Y}\sim\mathcal{DM}_{\Delta_N}(\boldsymbol{\alpha})\underset{N}{\wedge}\mathcal{L}(\psi)$ } \\
\hline
$\mathcal{L}(\psi)$
& $\mathcal{B}_{n}(p)$
& $\mathcal{NB}(r, p)$\footnote{If $r=\vert \boldsymbol{\alpha}\vert$ the splitting distribution is a not a \textit{sensu stricto} multivariate distribution because the graphical model is empty}
& $\mathcal{P}(\lambda)$
\\
\hline
Supp$(\boldsymbol{Y})$ 
& $\blacktriangle_{n}$ 
& $\mathbb{N}^J$
& $\mathbb{N}^J$
\\
\hline
$P(\boldsymbol{Y}=\boldsymbol{y})$ 
& $\frac{\Gamma(n+1)p^{\vert\boldsymbol{y}\vert}(1-p)^{n-\vert\boldsymbol{y}\vert}}{\Gamma(n-\vert\boldsymbol{y}\vert+1)(\vert\boldsymbol{\alpha}\vert)_{\vert\boldsymbol{y}\vert}}\prod_{j=1}^J\frac{(\alpha_j)_{y_j}}{y_j!}$
& $(1-p)^{r}\frac{(r)_{\vert\boldsymbol{y}\vert}p^{\vert\boldsymbol{y}\vert}}{(\vert\boldsymbol{\alpha}\vert)_{\vert\boldsymbol{y}\vert}}\prod_{j=1}^J\frac{(\alpha_j)_{y_j}}{y_j!}$
& $e^{-\lambda}\frac{\lambda^{\vert\boldsymbol{y}\vert}}{(\vert\boldsymbol{\alpha}\vert)_{\vert\boldsymbol{y}\vert}}\prod_{j=1}^J\frac{(\alpha_j)_{y_j}}{y_j!}$
\\
\hline
$\operatorname{E}(\boldsymbol{Y})$
& $\frac{np}{\vert\boldsymbol{\alpha}\vert}\cdot\boldsymbol{\alpha}$
& $\frac{rp}{\vert\boldsymbol{\alpha}\vert(1-p)}\cdot\boldsymbol{\alpha}$
& $\frac{\lambda}{\vert\boldsymbol{\alpha}\vert}\cdot\boldsymbol{\alpha}$
\\
\hline
$\operatorname{Cov}(\boldsymbol{Y})$
& 
\begin{tabular}{l} 
$\frac{np}{\vert\boldsymbol{\alpha}\vert(\vert\boldsymbol{\alpha}\vert+1)} \cdot\big\lbrace \left\lbrace (n-1)p + \vert\boldsymbol{\alpha}\vert +1 \right\rbrace \cdot\textrm{diag}(\boldsymbol{\alpha})$  \\ 
$ -  \frac{p(n+\vert\boldsymbol{\alpha}\vert)}{\vert\boldsymbol{\alpha}\vert}  \cdot\boldsymbol{\alpha}\boldsymbol{\alpha}^t \big\rbrace$ 
\end{tabular}
& \begin{tabular}{l} 
$\frac{rp}{\vert\boldsymbol{\alpha}\vert(\vert\boldsymbol{\alpha}\vert+1)(1-p)} \cdot\big[ \frac{(r-\vert\boldsymbol{\alpha}\vert)p+\vert\boldsymbol{\alpha}\vert+1}{1-p}  \cdot\textrm{diag}(\boldsymbol{\alpha})$  \\ 
$ + \frac{(\vert\boldsymbol{\alpha}\vert-r)p}{\vert\boldsymbol{\alpha}\vert(1-p)} \cdot\boldsymbol{\alpha}\boldsymbol{\alpha}^t \big]$ 
\end{tabular}
& $\frac{\lambda}{\vert\boldsymbol{\alpha}\vert(\vert\boldsymbol{\alpha}\vert+1)} \cdot \left\lbrace \left( \lambda + \vert\boldsymbol{\alpha}\vert +1 \right)\cdot\textrm{diag}(\boldsymbol{\alpha})  - \frac{\lambda}{\vert\boldsymbol{\alpha}\vert}  \cdot\boldsymbol{\alpha}\boldsymbol{\alpha}^t \right\rbrace$ 
\\
\hline
$G_{\boldsymbol{Y}}(\boldsymbol{s})$
& $(1-p)^{n}\;{}^J_1F_1(-n;\boldsymbol{\alpha};\vert\boldsymbol{\alpha}\vert;-\frac{p}{1-p}\cdot\boldsymbol{s})$
& $(1-p)^{r}\;{}^J_1F_1(r;\boldsymbol{\alpha};\vert\boldsymbol{\alpha}\vert;p\cdot\boldsymbol{s})$
& $e^{-\lambda}\;{}^J_0F_1(\boldsymbol{\alpha};\vert\boldsymbol{\alpha}\vert;\lambda\cdot\boldsymbol{s})$
\\
\hline
Marginals
& $Y_j\sim\beta_p\mathcal{B}_{n}(\alpha_j,\vert\boldsymbol{\alpha}_{-j}\vert)$
& $Y_j\sim\beta_p\mathcal{NB}(r,\alpha_j,\vert\boldsymbol{\alpha}_{-j}\vert)$
& $Y_j\sim\beta_{\lambda}\mathcal{P}(\alpha_j,\vert\boldsymbol{\alpha}_{-j}\vert)$
\\
\hline
\end{tabular}
\caption{\label{table:dirichlet:multinomial:psd}Usual characteristics 
of Dirichlet multinomial splitting binomial, negative binomial and Poisson distribution.}
\end{table}
\end{landscape}

\section{Splitting Regression models}
\label{section:regression}
Let us consider the regression framework, with the discrete multivariate response variable $\boldsymbol{Y}$ and the vector of $Q$ explanatory variables $\boldsymbol{X} = \left(X_1, \ldots , X_{Q}\right)$.
The random vector $\boldsymbol{Y}$ is said to follow a splitting regression if  there exists $\boldsymbol{\psi}: {\rm Supp}({\boldsymbol X}) \rightarrow {\boldsymbol \Psi}$ and $\boldsymbol{\theta}: {\rm Supp}({\boldsymbol X}) \rightarrow {\boldsymbol \Theta}$ such that:
\begin{itemize}
	\item the random vector $\boldsymbol{Y}$ given $\vert \boldsymbol{Y} \vert = n$ and $ \boldsymbol{X} = \boldsymbol{x}$ follows the singular regression $\mathcal{S}_{\Delta_{n}}\left\lbrace \boldsymbol{\theta}\left(\boldsymbol{x}\right)\right\rbrace$ for all $n \in \mathbb{N}$.
	\item the sum $\vert \boldsymbol{Y} \vert$ given $\boldsymbol{X} = \boldsymbol{x}$ follow the univariate regression $\mathcal{L}\left\lbrace \boldsymbol{\psi}\left(\boldsymbol{x}\right)\right\rbrace$.
\end{itemize}
Such a compound regression model will be denoted by $\left. \boldsymbol{Y} \,\middle\vert\, \boldsymbol{X} = \boldsymbol{x} \right. \sim \mathcal{S}_{\Delta_{N}}\left\lbrace \boldsymbol{\theta}\left(\boldsymbol{x}\right)\right\rbrace \underset{N}{\wedge} \mathcal{L}\left\lbrace\boldsymbol{\psi}\left(\boldsymbol{x}\right)\right\rbrace$.
The decomposition of log-likelihood ~\eqref{prop:inference:freq} still holds when considering covariates if parametrizations of the singular distribution and the sum distribution are unrelated.
Table~\ref{table:regressions} gives some references for parameter inference and variable selection adapted to three singular regressions and six univariate regressions. 
The choice of the link function for the singular regression is related to the symmetry of the resulting splitting regression. 
Using the singular multinomial regression for instance, only the canonical link function implies the symmetry of the splitting regression (see \citet{PTG15} for details about invariance properties of categorical regression models). 
Note that all generalized Dirichlet multinomial splitting regressions are not symmetric since the singular distribution is not. 

\begin{table}
(a) \\
   \begin{tabular}{|c|c|c|}
     \hline
Regression &  Link function & Parameter inference \\\hline\hline
Multinomial & 
$\pi_j=\frac{\exp(\boldsymbol{x}^t\boldsymbol{\beta}_j)}{1+\exp(\boldsymbol{x}^t\boldsymbol{\beta}_j)},\;j=1,\ldots,J-1$  & See \cite{ZZZS17} \\\hline
Dirichlet multinomial & $\alpha_j=\exp(\boldsymbol{x}^t\boldsymbol{\beta}_j),\;j=1,\ldots,J$ &
 see \cite{ZZZS17} \\
     	Generalized Dirichlet multinomial &        
        \begin{tabular}[t]{c}  $a_j=\exp(\boldsymbol{x}^t\boldsymbol{\beta}_{1,j}),\;j=1,\ldots,J-1$ \\  $b_j=\exp(\boldsymbol{x}^t\boldsymbol{\beta}_{2,j}),\;j=1,\ldots,J-1$ \end{tabular} & see \cite{ZZZS17}\\\hline
   \end{tabular}

(b) \\
   \begin{tabular}{|c|c|c|}
     \hline
    	Regression &  Link function & Parameter inference \\\hline\hline
    	Poisson & $\lambda=\exp(\boldsymbol{x}^t\boldsymbol{\beta})$ & See \cite{McCN89}  \\\hline
    	Binomial  & $p=\frac{\exp(\boldsymbol{x}^t\boldsymbol{\beta})}{1+\exp(\boldsymbol{x}^t\boldsymbol{\beta})}$ & See \cite{McCN89} for $n$ known\\\hline
    	Negative binomial & $p=\exp(\boldsymbol{x}^t\boldsymbol{\beta})$  & See \cite{Hil11} \\\hline\hline
        Beta Poisson & $\frac{a}{a+b}=\frac{\exp(\boldsymbol{x}^t\boldsymbol{\beta})}{1+\exp(\boldsymbol{x}^t\boldsymbol{\beta})}$  & See \cite{VWKNWRP16} \\\hline
     	Beta binomial & $\frac{a}{a+b}=\frac{\exp(\boldsymbol{x}^t\boldsymbol{\beta})}{1+\exp(\boldsymbol{x}^t\boldsymbol{\beta})}$ &
         \begin{tabular}[t]{c} See \cite{FF88} \\ and \cite{LL12} for $n$ known \end{tabular} \\
     	Beta negative binomial & $\frac{a}{a+b}=\frac{\exp(\boldsymbol{x}^t\boldsymbol{\beta})}{1+\exp(\boldsymbol{x}^t\boldsymbol{\beta})}$ &        
        \begin{tabular}[t]{c}  See \cite{RCSOM09} \\ and \cite{SVORCM17}  \end{tabular}
        \\\hline
   \end{tabular}
      \caption{\label{table:regressions} References of inference procedures for (a) singular regressions and (b) univariate regressions.}
\end{table}

\section{Empirical studies}
\label{section:application}

All studies presented in this section are reproducible.
Packages used are installable using the conda package management system and each study is available as a Jupyter notebook (see Appendix~\ref{CODE}).

\subsection{A comparison of multivariate models for count data}
\label{subsec:simemp:comparison}

In order to illustrate the variety of splitting models, we considered two datasets used in the literature to illustrate models for count data.
The first one consists in outcomes of football games \citep{KN05} and the second one consists in simulated data mimicking data obtained from sequencing techonologies such as RNA-seq data \citep{ZZZS17}.
The goal being to compare distributions and regressions models, comparisons were performed when considering all covariates or none of the covariates (see Table~\ref{table:comparison}).
Remark that variable selection \citep[e.g., using regularization methods]{ZZZS17} is possible, but is out of the scope of this paper.

Let us first remark that the inference methodology for multinomial, Dirichlet multinomial and generalized Dirichlet multinomial regressions presented by \cite{ZZZS17} and implemented by \cite{ZZ17} is only valid for singular versions.
Their comparisons of these models against the negative multinomial is therefore invalid since the first three models focus on $\boldsymbol{Y}$ given $\left\vert \boldsymbol{Y} \right\vert$ and the latter focuses on $\boldsymbol{Y}$.
Hence, we only compared our results to their unique $J$-multivariate model that is the negative multinomial model and the multivariate Poisson model defined by \cite{KM05}.
By limiting the number of sum models to $7$ and the number of singular models to $6$, we were able to propose $42$ splitting models.
Among those $42$ models, only $4$ models  were not \emph{sensu stricto} multivariate models since multinomial splitting Poisson models induce independent response variables. 

\begin{table}[t]
   \subfloat[]{\begin{tabular}{|c|c|c|}\hline
        $\boldsymbol{Y}$ given $\vert\boldsymbol{Y}\vert=n$ and $\boldsymbol{X} = \boldsymbol{x}$ 
        & BIC$_0$ & BIC$_1$\\\hline\hline
        $\mathcal{M}_{\Delta_n}\left(\boldsymbol{\pi}\right)$ & $574.18$ & $38,767.91$ \\\hline
        $\mathcal{DM}_{\Delta_n}\left(\boldsymbol{\alpha}\right)$ & $579.49$ & $9,969.121$ \\\hline
        $\mathcal{GDM}_{\Delta_n}\left(\boldsymbol{\alpha}, \boldsymbol{\beta}\right)$ &  $579.49$ & $9,735.45$ \\\hline
        $\mathcal{M}_{\Delta_n}\left\lbrace \boldsymbol{\pi}\left(\boldsymbol{x}\right)\right\rbrace$ & $\boldsymbol{508.14}$ & $15,145.24$ \\\hline
        $\mathcal{DM}_{\Delta_n}\left\lbrace \boldsymbol{\alpha}\left(\boldsymbol{x}\right)\right\rbrace$ & $836.4$ & $8,932.83$\\\hline
        $\mathcal{GDM}_{\Delta_n}\left\lbrace \boldsymbol{\alpha}\left(\boldsymbol{x}\right), \boldsymbol{\beta}\left(\boldsymbol{x}\right)\right\rbrace$ & $836.4$ & $\boldsymbol{8,843.479}$ \\\hline
   \end{tabular}}\quad
    \subfloat[]{\begin{tabular}{|c|c|c|}
        \hline
        $\vert\boldsymbol{Y}\vert$ given $\boldsymbol{X} = \boldsymbol{x}$ 
        & BIC$_0$ & BIC$_1$\\\hline\hline
        $\mathcal{P}\left(\lambda\right)$ & $\boldsymbol{1,130.64}$ & $13,074.12$ \\\hline
        $\mathcal{B}_n\left(p\right)$ & $1,165.6$ & $26,474.38$ \\\hline
        $\mathcal{NB}\left(r, p\right)$ & $1,131.85$ & $2,678.55$ \\\hline
        $\mathcal{L}\left(p\right)$ & $1,370.84$ & $3,513.92$ \\\hline
        $\mathcal{P}\left\lbrace\lambda\left(\boldsymbol{x}\right)\right\rbrace$ & $1,258.65$ & $6,353.13$ \\\hline
        $\mathcal{B}_n\left\lbrace p\left(\boldsymbol{x}\right)\right\rbrace$ & $1,272.7$ & $12,999.47$ \\\hline
        $\mathcal{NB}\left\lbrace r, p\left(\boldsymbol{x}\right)\right\rbrace$ & $1,264.38$ & $\boldsymbol{2,514.30}$ \\\hline
   \end{tabular}}\\
   \subfloat[]{\begin{tabular}{|c|c|c|}
     \hline
        $\boldsymbol{Y}$ given $\boldsymbol{X} = \boldsymbol{x}$ & BIC$_0$ & BIC$_1$\\\hline\hline
        $\mathcal{MP}\left(\boldsymbol{\lambda}\right)$ &  $\boldsymbol{1,710.05}$ & $\diagup\!\!\!\!\!\diagdown$ \\
        $\mathcal{MP}\left\lbrace\boldsymbol{\lambda}\left(\boldsymbol{x}\right)\right\rbrace$ &  $1,956.10$ & $\diagup\!\!\!\!\!\diagdown$ \\\hline
   \end{tabular}}\quad
   \subfloat[]{\begin{tabular}{|c|c|c|}
     \hline
        $\boldsymbol{Y}$ given $\boldsymbol{X} = \boldsymbol{x}$ & BIC$_0$ & BIC$_1$\\\hline\hline
        $\mathcal{MN}\left(r, \boldsymbol{\pi}\right)$ & $\boldsymbol{1,705.93}$ & $41,384.52$\\\hline
        $\mathcal{MN}\left\lbrace r, \boldsymbol{\pi}\left(\boldsymbol{x}\right)\right\rbrace$ & $2,176.3$ & $\boldsymbol{17,657.63}$ \\\hline
   \end{tabular}}
   \caption{\label{table:comparison}Bayesian Information Criteria (BIC) obtained for the first dataset \citep[BIC$_0$]{KN05} and the second one \citep[BIC$_1$]{ZZZS17} for (a) singular models, (b) sum models , (c) Poisson and (d) negative multinomial models.
           Multivariate Poisson models could not be fit to the second dataset since, to our knowledge, there is no implementation available in R for more than $2$ response variables \citep{KN05}.}
   \label{table}
\end{table}

For the first dataset, the best splitting model is a singular multinomial regression compounded by a Poisson distribution with a BIC of $508.14+1,130.64=1,638.78$.
This score is inferior to the one of the best multivariate Poisson model (i.e., $1,710.05$) and the one of the best negative multinomial model (i.e., $1,705.93$).
This indicates that there is no relationship between football team goals.
For the second dataset, the best splitting model is a singular generalized Dirichlet multinomial regression compounded by a negative binomial regression with a BIC of $8843.48+2514.3=11,357.78$.
This score is also inferior to the one of the best negative multinomial model (i.e., $17,657.63$).

\subsection{An application to mango patchiness analysis}
\label{subsec:simemp:application}

Recently, a statistical methodology has been proposed to characterize plant patchiness at the plant scale \citep{FDDPLNG16}.
However, little is known about patchiness at the whole population scale.
To characterize patchiness at the plant scale, a segmentation/clustering of tree-indexed data method has been proposed in order to split an heterogeneous tree into multiple homogeneous subtrees.
After the clustering, 
the tree can be summarized into a multivariate count denoting the number of subtrees in each cluster (i.e., patch type).
Mixture of multinomial splitting distributions can therefore be considered to recover the different types of tree patchiness that can be found in the plant population.
Such a mixture model is of high interest since it enables to discriminate the types of tree patchiness according to the :
\begin{itemize}
	\item number of patches present on trees, by fitting different sum distributions within components of the mixture model,
    \item distribution of these patches among types, by fitting different singular distributions within components of the mixture model.
\end{itemize}
We here consider results presented by \cite{FDDPLNG16} to conduct our \textit{post-hoc} analysis.
Three different types of patches have been identified for mango trees: vegetative patches which contain almost only vegetative growth units (GU, plant elementary component), reproductive patches which contain almost only GUs that flowered or fructified and quiescent patches which contain GUs that did not burst, flowered nor fructified.
Multinomial splitting distributions of mixture components were therefore of dimension $3$, where $N_0$ (resp. $N_1$ and $N_2$) denotes the number of vegetative (resp. reproductive and quiescent) patches observed within a tree.
Since there is at least one patch in a mango tree (i.e., the tree itself), shifted singular multinomial splitting distributions were considered with a shift equal to $1$ for binomial, negative binomial and Poisson sum distributions but without shift for geometric and logarithmic distributions.
Within each component the parametric form of the sum distribution was selected using the BIC.

\begin{figure}
	\begin{center}
     \scalebox{.6}{\input{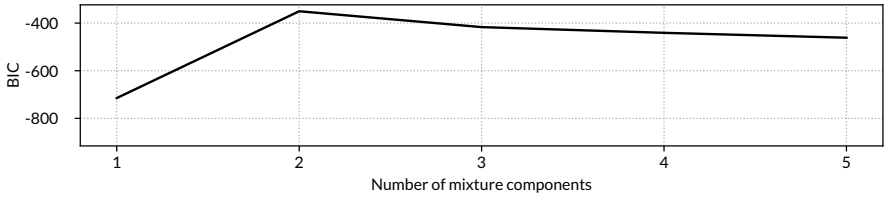}}
    \caption{BIC according to the number of multinomial splitting components of mixtures.}
    \label{fig:application:bic}
    \end{center}
\end{figure}

\begin{figure}
	\begin{center}
    \scalebox{.8}{\input{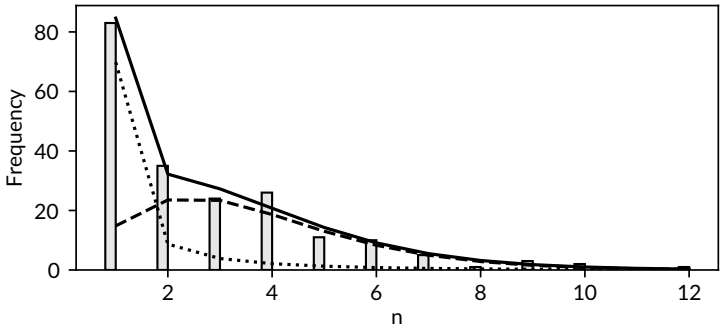}}\label{fig:application:fit}
    \caption{Representation of the mixture of sum distributions estimated (with a solid line) confronted to data frequencies (gray bars).
    Note that the sum distribution of the first (resp. second) component is represented with a dotted (resp. dashed) line.}
    \label{fig:application:mixture}
    \end{center}
\end{figure}

The mixture model selected using BIC has two components (see Figure~\ref{fig:application:bic}) with weights $P\left(L = 1\right) = 0.44$ and $P\left(L = 2\right) = 0.56$.
In the two components $i=1,2$, the number of patches followed a multinomial splitting shifted negative binomial distribution $\boldsymbol{Y} \,\vert\, L = i \sim \mathcal{M}_{\Delta_N}\left(\boldsymbol{\pi}_i\right) \underset{N}{\wedge} \mathcal{NB}\left(r_i, p_i; \delta_i\right)$ with estimations $\hat{\boldsymbol{\pi}}_1=\left(0.21, 0.00, 0.79\right)$, $\hat{r}_1=0.16$, $\hat{p}_1=0.76$, $\hat{\delta}_1=1$ for the first component and $\hat{\boldsymbol{\pi}}_2=\left(0.54, 0.17, 0.28\right)$, $\hat{r}_2=3.96$, $\hat{p}_2=0.40$, $\hat{\delta}_2=1$ for the second component.
This mixture of two components indicates that the population of mango trees can be separated into two types of trees (see Figure~\ref{fig:application:mixture}):
\begin{itemize}
    \item mango trees with a relatively low number of patches that can be either vegetative or quiescent but not reproductive (component 1),
	\item mango trees with a relatively high number of patches that can be of any type and in particular reproductives (component 2).
\end{itemize}
These types of trees are almost equally represented in the period considered ($52\%$ for the first component against $48\%$).
This result tends to imply that the reproductive period of mango trees leads to an increase in patch number increase while the vegetative period leads to a decrease in patch number.


\section{Discussion}
\label{sec:discussion}

Convolutions splitting distributions that are positive and additive, have been studied in depth in this paper since their graphical models and their marginal distributions are easily obtained.
The characterization of the graphical model of hypergeometric splitting distributions stay an open issue because of the non-positivity.
But thanks to the additivity, Theorem \ref{marginal_csd} still holds.
It would be interesting to instantiate some univariate distributions $\mathcal{L}(\psi)$ and precisely describe the resulting splitting distributions.
More generally, the multivariate Polya distribution with parameters $n\in\mathbb{N}$, $\boldsymbol{\theta}\in\boldsymbol{\Theta}$ and $c\in\mathbb{R}$ encompasses the multivariate hypergeometric ($c=-1$), the multinomial ($c=0$) and the Dirichlet multinomial ($c=1$) distributions \citep{JP70}.
It would therefore be interesting to study the properties of multivariate Polya splitting distributions according to the $c$ value.
Otherwise, non-symmetric convolution distributions could be defined (including the generalized Dirichlet multinomial distribution as a special case) to ease the study of corresponding splitting distributions.

Another alternative to define new singular distributions is to consider their mixture.
Finite mixture can be inferred using classical expectation-maximization algorithm for multivariate distributions.
Moreover, in an application context the principle of mixture models is quite interesting for splitting models.
If we consider the mango tree application, we inferred mixture of splitting distributions in order to characterize plant patchiness at the plant scale.
This relied on the assumption that tree patchiness is both expressed in terms of number of patches and their type distribution.
One the one hand, if tree patchiness is only a phenomenon expressed in term of number of patches, a mixture of sum distributions could be considered to distinguish trees.
On the other hand, if tree patchiness is only a phenomenon expressed in term of patch type distribution, singular distributions constructed using mixture of singular distributions could be of most interest.

Finally, this work could be used for learning graphical models with discrete variables, which is an open issue.
Although the graphical models for convolution splitting distributions are basic (complete or empty), they could be used as building blocks for partially directed acyclic graphical models. 
Therefore, the procedure of learning partially directed acyclic graphical models described by \citet{FDG14} could be used for learning graphical models based on convolution splitting distributions and regressions.
It could be used for instance to infer gene co-expression network from RNA seq dataset. 

\bibliographystyle{elsarticle-harv}
\bibliography{references.bib}


\appendix
\section{\label{proofs} Proofs}

\paragraph{Details about marginal of symmetric splitting distributions}
For any $j\in\{1,\ldots,J\}$ and $y_j \in \mathbb{N}$ we have
\begin{align*}
P(Y_j=y_j) & = \sum_{\boldsymbol{y}_{-j}} P(\boldsymbol{Y} =\boldsymbol{y}),\\
& =   \sum_{n \geq y_j }   P(\vert \boldsymbol{Y} \vert = n)      \sum_{\boldsymbol{y}_{-j}} P_{\vert \boldsymbol{Y} \vert = n}(\boldsymbol{Y} =\boldsymbol{y} ),   \\
P(Y_j=y_j)  & =   \sum_{n \geq y_j }   P(\vert \boldsymbol{Y} \vert = n)  P_{\vert \boldsymbol{Y} \vert = n}(Y_j =y_j), 
\end{align*}
The marginal distribution of the singular distribution, i.e., the distribution of $Y_j$ given $\vert \boldsymbol{Y} \vert = n$, is a distribution bounded by $n$.
Its parametrization has the same form $f_j(\boldsymbol{\theta})$ for all marginals $Y_j$ given $\vert \boldsymbol{Y} \vert = n$, since the singular distribution is symmetric.
It implies that all marginals $Y_j$ follow the damage distribution $\mathcal{L}_N \{f_j(\boldsymbol{\theta})\}\underset{N}{\wedge} \mathcal{L}(\psi)$.

\paragraph{Proof of Theorem \ref{marginal_csd}}
Let $\mathcal{I}$ denote a subset of $\{1,\ldots,J\}$ with$\mathcal{I}\neq\{1,\ldots,J\}$.
Let $-\mathcal{I}$ denote the subset $\{1,\ldots,J\}\setminus\mathcal{I}$ and $\boldsymbol{y}_{\mathcal{I}}$ (respectively $\boldsymbol{y}_{-\mathcal{I}}$) denote the corresponding sub-vectors.

\subparagraph{Proof of \ref{item:marginal:sum}}  
\begin{align*}
P(\vert\boldsymbol{Y}_{\mathcal{I}}\vert=n) & =  \sum_{\boldsymbol{y}_{\mathcal{I}}\in\Delta_n} P(\boldsymbol{Y}_{\mathcal{I}}=\boldsymbol{y}_{\mathcal{I}})\\
& = \sum_{\boldsymbol{y}_{\mathcal{I}}\in\Delta_n} \sum_{\boldsymbol{y}_{-\mathcal{I}}} P(\boldsymbol{Y}=\boldsymbol{y})\\
& =\sum_{k\geq n} P(\vert\boldsymbol{Y}\vert=k) \sum_{\boldsymbol{y}_{\mathcal{I}}\in\Delta_n} \sum_{\boldsymbol{y}_{-\mathcal{I}}\in\Delta_{k-n}} P_{\vert\boldsymbol{Y}\vert=k}(\boldsymbol{Y}=\boldsymbol{y})  \\
& = \sum_{k\geq n}    \frac{P(\vert\boldsymbol{Y}\vert=k)}{c_{\boldsymbol{\theta}}(k)} \sum_{\boldsymbol{y}_{\mathcal{I}}\in\Delta_n} \prod_{j\in\mathcal{I}} a_{\theta_j}(y_j) \sum_{\boldsymbol{y}_{-\mathcal{I}}\in\Delta_{k-n}}   \prod_{j\in -\mathcal{I}} a_{\theta_j}(y_j) \\
P(\vert\boldsymbol{Y}_{\mathcal{I}}\vert=n) & = \sum_{k\geq n}  \frac{ c_{\boldsymbol{\theta}_{\mathcal{I}}} (n) c_{\boldsymbol{\theta}_{-\mathcal{I}}}(k-n)}{c_{\boldsymbol{\theta}}(k)} P(\vert\boldsymbol{Y}\vert=k)
\end{align*}
where $c_{\boldsymbol{\theta}_\mathcal{I}}(n)$ denotes the convolution of $(a_{\theta_j})_{j\in \mathcal{I}}$ over the simplex $\Delta_n$. 
Since the convolution distribution is assumed to be additive, we obtain by recursion on $j\in \mathcal{I}$ (resp. $j\in -\mathcal{I}$ and $j\in\{1,\ldots,J\}$) that
\begin{equation}\label{marginal:sum:distribution}
P(\vert\boldsymbol{Y}_{\mathcal{I}}\vert=n) = \sum_{k\geq n}  \frac{ a_{\vert\boldsymbol{\theta}_{\mathcal{I}}\vert} (n) a_{\vert\boldsymbol{\theta}_{-\mathcal{I}}\vert}(k-n)}{a_{\vert\boldsymbol{\theta}\vert}(k)} P(\vert\boldsymbol{Y}\vert=k)
\end{equation}
Moreover we obtain the convolution identity $\sum_{n=0}^{k} a_{\vert\boldsymbol{\theta}_\mathcal{I}\vert}(n) a_{\vert\boldsymbol{\theta}_{-\mathcal{I}}\vert}(k-n) =a_{\vert\boldsymbol{\theta}\vert}(k)$ and thus the last equation well defines the desired convolution damage distribution  $\sim\mathcal{C}_{N}(a; \vert \boldsymbol{\theta}_{\mathcal{I}} \vert, \vert \boldsymbol{\theta}_{-\mathcal{I}} \vert ) \underset{N}{\wedge} \mathcal{L}(\boldsymbol{\psi})$.

\subparagraph{Proof of \ref{item:marginal:singular} and \ref{item:marginal}}
\begin{align*}
P(\boldsymbol{Y}_{\mathcal{I}}=\boldsymbol{y}_{\mathcal{I}}, \; \vert\boldsymbol{Y}_{\mathcal{I}}\vert=n) & = P(\boldsymbol{Y}_{\mathcal{I}}=\boldsymbol{y}_{\mathcal{I}}) \mathbb{1}_{\Delta_n}(\boldsymbol{y}_{\mathcal{I}}) \\
& = \mathbb{1}_{\Delta_n}(\boldsymbol{y}_{\mathcal{I}}) \sum_{k\geq n} P(\vert\boldsymbol{Y}\vert=k)  \sum_{\boldsymbol{y}_{-\mathcal{I}}\in\Delta_{k-n}} P_{\vert\boldsymbol{Y}\vert=k}(\boldsymbol{Y}=\boldsymbol{y})   \\
& = \prod_{j\in\mathcal{I}} a_{\theta_j}(y_j) \mathbb{1}_{\Delta_n}(\boldsymbol{y}_{\mathcal{I}}) \sum_{k\geq n}    \frac{P(\vert\boldsymbol{Y}\vert=k)}{c_{\boldsymbol{\theta}}(k)}  \sum_{\boldsymbol{y}_{-\mathcal{I}}\in\Delta_{k-n}}   \prod_{j\in -\mathcal{I}} a_{\theta_j}(y_j) \\
P(\boldsymbol{Y}_{\mathcal{I}}=\boldsymbol{y}_{\mathcal{I}}, \; \vert\boldsymbol{Y}_{\mathcal{I}}\vert=n) & =  \prod_{j\in\mathcal{I}} a_{\theta_j}(y_j) \mathbb{1}_{\Delta_n}(\boldsymbol{y}_{\mathcal{I}})  \sum_{k\geq n}  \frac{a_{\vert\boldsymbol{\theta}_{-\mathcal{I}}\vert}(k-n)}{a_{\vert\boldsymbol{\theta}\vert}(k)} P(\vert\boldsymbol{Y}\vert=k)
\end{align*}
Using equation \ref{marginal:sum:distribution} we obtain the conditional probability 
\[
P_{\vert\boldsymbol{Y}_{\mathcal{I}}\vert=n}(\boldsymbol{Y}_{\mathcal{I}}=\boldsymbol{y}_{\mathcal{I}}) = \frac{1}{a_{\vert\boldsymbol{\theta}_{\mathcal{I}}\vert} (n)} \prod_{j\in\mathcal{I}} a_{\theta_j}(y_j) \mathbb{1}_{\Delta_n}(\boldsymbol{y}_{\mathcal{I}}),
\]
and thus \ref{item:marginal:singular} holds.
Remark that \ref{item:marginal:sum} and \ref{item:marginal:singular} imply \ref{item:marginal} by definition of a splitting distribution.

\subparagraph{Proof of \ref{item:full:conditional} and \ref{item:general:conditional}}
\begin{align*}
P_{\boldsymbol{Y}_{-\mathcal{I}}=\boldsymbol{y}_{-\mathcal{I}}}(\boldsymbol{Y}_{\mathcal{I}}=\boldsymbol{y}_{\mathcal{I}}) & =  P_{\boldsymbol{Y}_{-\mathcal{I}}=\boldsymbol{y}_{-\mathcal{I}},\; \vert\boldsymbol{Y}_{\mathcal{I}}\vert= \vert\boldsymbol{y}_{\mathcal{I}}\vert}(\boldsymbol{Y}_{\mathcal{I}}=\boldsymbol{y}_{\mathcal{I}}) P_{\boldsymbol{Y}_{-\mathcal{I}}=\boldsymbol{y}_{-\mathcal{I}}}(\vert\boldsymbol{Y}_{\mathcal{I}}\vert= \vert\boldsymbol{y}_{\mathcal{I}}\vert) 
\end{align*}
Since the sum $\vert\boldsymbol{Y}\vert$ is independent of the vector $\boldsymbol{Y}_{-\mathcal{I}}$ given its sum $\vert\boldsymbol{Y}_{-\mathcal{I}}\vert$ it can be shown that
\begin{align*}
P_{\boldsymbol{Y}_{-\mathcal{I}}=\boldsymbol{y}_{-\mathcal{I}}}(\boldsymbol{Y}_{\mathcal{I}}=\boldsymbol{y}_{\mathcal{I}}) & =  P_{\vert\boldsymbol{Y}_{\mathcal{I}}\vert= \vert\boldsymbol{y}_{\mathcal{I}}\vert}(\boldsymbol{Y}_{\mathcal{I}}=\boldsymbol{y}_{\mathcal{I}}) 
 P_{\vert \boldsymbol{Y}_{-\mathcal{I}}\vert =\vert\boldsymbol{y}_{-\mathcal{I}}\vert}(\vert\boldsymbol{Y}_{\mathcal{I}}\vert= \vert\boldsymbol{y}_{\mathcal{I}}\vert) 
\end{align*}
Thanks to the result \ref{item:marginal:singular}, the left part of this product is given by the singular convolution distribution.
Remarking that $P_{\vert \boldsymbol{Y}_{-\mathcal{I}}\vert =\vert\boldsymbol{y}_{-\mathcal{I}}\vert}(\vert\boldsymbol{Y}_{\mathcal{I}}\vert= \vert\boldsymbol{y}_{\mathcal{I}}\vert) = P_{\vert\boldsymbol{Y}\vert\geq a}(\vert\boldsymbol{Y}\vert= a+ \vert\boldsymbol{y}_{\mathcal{I}}\vert) $ with $a=\vert\boldsymbol{y}_{-\mathcal{I}}\vert$ the left part is given by the truncated and shifted distribution $TS_a\{\mathcal{L}(\psi)\}$ and thus \ref{item:full:conditional} holds. 
Remark that \ref{item:marginal} and \ref{item:full:conditional} imply \ref{item:general:conditional}.

\paragraph{Proof of corollary \ref{coro:natural:extension}}
Assume that $\mathcal{L}(\boldsymbol{\psi})$ is stable under the damage process $\mathcal{C}_{N}(a; \vert \boldsymbol{\theta}_{\mathcal{I}} \vert, \vert \boldsymbol{\theta}_{-\mathcal{I}} \vert ) \underset{N}{\wedge} (\cdot)$ for any subset $\mathcal{I}\subset\{1,\ldots,J\}$.
Thanks to the additivity of the convolution distribution, Theorem \ref{marginal_csd} can be applied.
Using the property \ref{item:marginal}, it is easily seen that multivariate marginals are stable.
The criterion \ref{criter4} holds and the convolution splitting distribution is considered as a natural multivariate extension of $\mathcal{L}(\boldsymbol{\psi})$.
In particular, $\mathcal{L}(\boldsymbol{\psi})$ is stable under $\mathcal{C}_{N}(a; \vert \theta_{j} \vert, \vert \boldsymbol{\theta}_{-j} \vert ) \underset{N}{\wedge} (\cdot)$ and thus the univariate marginal follow $\mathcal{L}(\boldsymbol{\psi}_j)$ for some $\boldsymbol{\psi}_j\in\boldsymbol{\Psi}$.

\paragraph{Proof of corollary \ref{coro:convolution:nonsingular}}
Let $\boldsymbol{Y}$ follow the non-singular version of an additive convolution distribution: $\boldsymbol{Y}\sim\mathcal{C}_{\blacktriangle_n}(a_{\theta};\boldsymbol{\theta},\theta)$.
It means that the completed vector $(\boldsymbol{Y}, n-\vert\boldsymbol{Y}\vert)$ follow the additive convolution $\mathcal{C}_{\Delta_n^{J+1}}(a_{\theta};\boldsymbol{\theta},\theta)$.
Otherwise this singular distribution can seen as a particular splitting Dirac distribution, i.e.,  $\mathcal{C}_{\Delta_n^{J+1}}(a_{\theta};\boldsymbol{\theta},\theta) = \mathcal{C}_{\Delta_N^{J+1}}(a_{\theta};\boldsymbol{\theta},\theta) \underset{N}{\wedge} \mathbb{1}_n$.
Thanks to the additivity, the Theorem \ref{marginal_csd} can be applied on the completed vector $(\boldsymbol{Y}, n-\vert\boldsymbol{Y}\vert)$ to describe the distribution of $\boldsymbol{Y}$ (property \ref{item:marginal}):
\begin{align*}
    	\boldsymbol{Y} & \sim \mathcal{C}_{\Delta_{N}}\left(a_{\theta}; \boldsymbol{\theta}\right)  \underset{N}{\wedge}  \left\lbrace \mathcal{C}_{N^{\prime}}(a_{\theta}; \vert \boldsymbol{\theta} \vert, \theta ) \underset{N^{\prime}}{\wedge} \mathbb{1}_n \right\rbrace, \\
     \Leftrightarrow   \boldsymbol{Y} & \sim \mathcal{C}_{\Delta_{N}}\left(a_{\theta}; \boldsymbol{\theta}\right)  \underset{N}{\wedge}  \mathcal{C}_{n}(a_{\theta}; \vert \boldsymbol{\theta} \vert, \theta).
\end{align*}

\paragraph{Proof of equality \eqref{msd}}
Using the Corollary \ref{coro:convolution:nonsingular} with $a_{\theta}(y)=\theta^y/y!$ we obtain for $\boldsymbol{\theta}\in(0,\infty)^J$ and $\gamma\in(0,\infty)$
\begin{equation*}
\mathcal{M}_{\Delta_N}\left( \boldsymbol{\theta}\right) \underset{N}{\wedge} \mathcal{B}_{n}\left( \vert\boldsymbol{\theta}\vert,  \gamma\right)  = \mathcal{M}_{\blacktriangle_{n}}\left( \boldsymbol{\theta},  \gamma\right).
\end{equation*}
Denoting by $\boldsymbol{\pi}=\frac{1}{\vert\boldsymbol{\theta}\vert} \cdot\boldsymbol{\theta}$, $p=\frac{\vert\boldsymbol{\theta}\vert}{\vert\boldsymbol{\theta}\vert+\gamma}$ and $\boldsymbol{\pi}^*=\frac{1}{\vert\boldsymbol{\theta}\vert+\gamma} \cdot\boldsymbol{\theta}$ and using the proportionality we obtain equivalently
\begin{equation*}
\mathcal{M}_{\Delta_N}\left(\boldsymbol{\pi}\right) \underset{N}{\wedge} \mathcal{B}_{n}\left(p,1-p\right)  = \mathcal{M}_{\blacktriangle_{n}}\left( \boldsymbol{\pi}^*, 1-\vert \boldsymbol{\pi}^*\vert \right).
\end{equation*}
The notation of the binomial and the non-singular multinomial are then simplified by letting aside the last parameter without loss of generality, i.e. we have $\mathcal{M}_{\Delta_N}\left( \boldsymbol{\pi}\right) \underset{N}{\wedge} \mathcal{B}_{n}\left( p\right) = \mathcal{M}_{\blacktriangle_{n}}\left(\boldsymbol{\pi}^*\right)$.
Finally remarking that $\boldsymbol{\pi}^*=p\cdot\boldsymbol{\pi}$ we obtain the desired result.

\section{Remarkable continuous univariate distribution}
\label{RCUD}
Let us recall the definition of the (generalized) beta distribution with positive real parameters $c$, $\alpha$ and $b$, denoted by $\beta_{c}\left(a, b\right)$.
Its probability density function described by \citet{Whi72} is given by
\begin{equation*}
      f\left(x\right) = \frac{1}{B\left(a, b\right)} \frac{x^{a-1} \left(c-x\right)^{b-1}}{c^{a + b - 1}} \cdot \mathbb{1}_{\left(0, c\right)}(x).
\end{equation*}
Note that $Z = d X$ with $d \in (0,\infty)$  and $X \sim \beta_{c}\left(a, b\right)$ implies that $Z \sim \beta_{cd}\left(a, b\right)$.
The parameter $c$ of the beta distribution can thus be interpreted as a rescaling parameter of the standard beta distribution.
By convention the standard beta distribution (i.e., defined with $c=1$) will be denoted by $\beta \left(a, b\right)$.

Let us introduce the definition of the (generalized) beta square distribution with parameters $(a_1,b_1,a_2,b_2)\in(0,\infty)^4$ and $c\in(0,\infty)$, denoted by $\beta^2_{c}(a_1,b_1,a_2,b_2)$, as the product of the two independent beta distributions $\beta(a_1,b_1)$ and $\beta(a_2,b_2)$ normalized on $(0,c)$; see \citet{Dun13} for details.
It is named the standard beta square distribution when $c=1$ and denoted by $\beta^2(a_1,b_1,a_2,b_2)$. 
More generally the product of $m$ beta distributions could be defined.

\section{Remarkable discrete univariate distributions}
\label{RDUD}

\subsection{Power series distributions}
Let $(b_y)_{y\in\mathbb{N}}$ be a non-negative real sequence such that the series $\sum_{y \geq 0} b_y \theta^y$ converges toward $g(\theta)$ for all $\theta \in D=(0,R)$, where $R$ is the radius of convergence. The discrete random variable $Z$ is said to follow a power series distribution if for all $y \in \mathbb{N}$
\begin{equation*} P(Y=y) = \frac{b_y \theta^y}{g(\theta)}, \end{equation*}
and is denoted by $Y \sim PSD\{g(\theta)\}$.
Several usual discrete distributions fall into the family of power series distributions:
\begin{enumerate}
\item The Poisson distribution $\mathcal{P}(\lambda)$ with $b_y=1/y!$, $\theta=\lambda$, $g(\theta)=e^{\theta}$ and $D=(0,\infty)$.
\item The binomial distribution $\mathcal{B}_{n}(p)$ with $b_y= \binom{n}{y} \boldsymbol{1}_{y\leq n}$, $\theta=p/(1-p)$, $g(\theta)=(1+\theta)^{n}$ and $D=(0,\infty)$.
\item The negative binomial distribution $\mathcal{NB}(r, p)$ with $b_y= \binom{r + y-1}{y}$, $\theta=p$, $g(\theta)=(1-\theta)^{-r}$ and $D=(0,1)$.
\item The geometric distribution $\mathcal{G}(p)$  with $b_y= \boldsymbol{1}_{y\geq 1}$, $\theta=1-p$, $g(\theta)=\theta/(1-\theta)$ and $D=(0,1)$.
\item The logarithmic series distribution $\mathcal{L}(p)$  with $b_y= \boldsymbol{1}_{y\geq 1} 1/y$, $\theta=p$, $g(\theta)=- \ln (1-\theta)$ and $D=(0,1)$.
\end{enumerate}
When the support is a subset of $\mathbb{N}$, the $b_y$ values can be weighted by an indicator function as for binomial, geometric and logarithmic distributions. The $b_y$ must be independent of $\theta$ but they may depend on other parameters as for binomial and negative binomial distributions.

\subparagraph{Zero modified logarithmic series}\cite{JKB97}
The discrete variable $Y$ is said to follow a zero modified logarithmic series distribution with parameter $\omega\in [0,1)$ and $p \in (0,1)$ if the its probabilities are given by
\begin{align*}
P(Y=0) & = \omega, \\
P(Y=y) & = \frac{(1-\omega)p^y}{-y \ln (1-p)}, \; y \geq 1.
\end{align*}
This distribution will be denoted by $\mathcal{L}(p,\omega)$. Note that if $\omega=0$ this is the logarithmic series distribution: $\mathcal{L}(p,0)=\mathcal{L}(p)$.

\subsection{Beta compound distributions}
Usual characteristics of the standard beta binomial \citep{TGRGJ94}, standard beta negative binomial - also described by \citet{Xek81} as the univariate generalized waring distribution (UGWD) -  and the beta Poisson distributions \citep{Gur58} are first recalled in Table \ref{table:standard:beta:compound}.
Then we introduce these beta compound distributions in a general way, i.e. using the generalized beta distribution (see Appendix \ref{RCUD}).
For the Poisson case we obtain the same distribution since $\mathcal{P}(\lambda p) \underset{p}{\wedge} \beta(a,b) = \mathcal{P}(\theta) \underset{\theta}{\wedge} \beta_{\lambda}(a,b)$.
The two other case lead us to new distributions (\ref{table:generalized:beta:compound} for usual characteristics).
Remark that if $\pi=1$ then, the generalized beta binomial (resp. generalized beta negative binomial) turns out to be the standard beta binomial (resp. standard negative binomial distribution).
In opposition, if $\pi<1$, the non-standard beta binomial distribution (respectively non-standard beta negative binomial distribution) is obtained.

\paragraph{Generalized beta binomial distribution}
Let $n\in\mathbb{N}$, $a\in(0,\infty)$, $b\in(0,\infty)$ and $\pi\in (0,1)$ and consider the compound distribution $\mathcal{B}_{n}\left( p\right) \underset{p}{\wedge} \beta_{\pi}\left(a, b\right)$ denoted by $\beta_{\pi}\mathcal{B}_{n}\left(a, b\right)$.
Considering $\pi$ as a rescaling parameter, we have $\beta_{\pi}\mathcal{B}_{n}\left(a, b\right) = \mathcal{B}_{n}\left(\pi p\right) \underset{p}{\wedge} \beta\left(a, b\right)$.
Moreover, using the \gls{pgf} of the binomial distribution in equation \eqref{binomial-operator}, it can be shown that $\mathcal{B}_{n}(\pi p)=\mathcal{B}_{N}(\pi) \underset{N}{\wedge} \mathcal{B}_{n}(p)$.
Finally, using the Fubini theorem we obtain 
\begin{align*}
\beta_{\pi}\mathcal{B}_{n}(a,b) & =\left\lbrace \mathcal{B}_{N}(\pi) \underset{N}{\wedge} \mathcal{B}_{n}(p) \right\rbrace \underset{p}{\wedge} \beta(a,b), \\
& = \mathcal{B}_{N}(\pi) \underset{N}{\wedge} \left\lbrace \mathcal{B}_{n}(p) \underset{p}{\wedge} \beta(a,b) \right\rbrace, \\
\beta_{\pi}\mathcal{B}_{n}(a,b) & = \mathcal{B}_{N}(\pi) \underset{N}{\wedge} \beta\mathcal{B}_{n}\left(a, b\right).
\end{align*}
This is a binomial damage distribution whose the latent variable $N$ follows a standard beta binomial distribution.
The equation \eqref{binomial-operator} can thus be used to compute the probability mass function.
The $y^{th}$ derivative of the \gls{pgf} of the standard beta binomial distribution is thus needed
\[
	G_N^{(y)}(s) = \frac{(b)_n}{(a+b)_n} \frac{(-n)_y(a)_y}{(-b-n+1)_y} {}_2F_1\{(-n+y, a+y); -b-n+1+y; s\},
\]
obtained by induction on $y \in \mathbb{N}$.
The moments are obtained with the total law of expectation and variance given the latent variable $N$ of the binomial damage distribution.
In the same way, we obtain the \gls{pgf} as $G_Y(s)=G_N(1-\pi+\pi s)$.
A similar proof holds for the generalized beta negative binomial case.

\paragraph{Generalized beta square compound distributions}
It is also possible to define the (generalized) beta square distribution, as the product of two independent beta distributions \citep{Dun13}, and then define the (generalized) beta square compound distributions.

\begin{itemize}
\item The standard beta square binomial distribution is defined as $\mathcal{B}_n(p)\underset{p}{\wedge}\beta^2(a_1,b_1,a_2,b_2)$ and denoted by $\beta^2\mathcal{B}_n(a_1,b_1,a_2,b_2)$.
\item The standard beta square negative binomial distribution is defined as $\mathcal{NB}(r,p)\underset{p}{\wedge}\beta^2(a_1,b_1,a_2,b_2)$ and denoted by $\beta^2\mathcal{NB}(r,a_1,b_1,a_2,b_2)$.
\item The generalized beta square binomial distribution is defined as $\mathcal{B}_n(p)\underset{p}{\wedge}\beta^2_{\pi}(a_1,b_1,a_2,b_2)$ and denoted by $\beta^2_{\pi}\mathcal{B}_n(a_1,b_1,a_2,b_2)$ 
\item The generalized beta square negative binomial distribution is defined as $\mathcal{NB}(r,p)\underset{p}{\wedge}\beta^2_{\pi}(a_1,b_1,a_2,b_2)$ and denoted by $\beta^2_{\pi}\mathcal{NB}(r,a_1,b_1,a_2,b_2)$ 
\item The beta square Poisson distribution is defined as $\mathcal{P}(\theta)\underset{\theta}{\wedge}\beta_{\lambda}^2(a_1,b_1,a_2,b_2)$ and denoted by $\beta_{\lambda}^2\mathcal{P}(a_1,b_1,a_2,b_2)$.
\end{itemize}

\section{Reproducibility}
\label{CODE}

The source code (written in {\it C++} and {\it Python}) used for the inference of splitting distributions is available on GitHub (\url{https://github.com/StatisKit}) and binaries can be installed using the Conda package management system (\url{(http://conda.pydata.org)}).
Refers to the documentation for more information (\url{(http://statiskit.rtfd.io)}).

Our analyses performed with the {\it Python} interface  or {\it R} packages is available in the Jupyter notebook format as supplementary materials and can be reproduced using the Docker \citep{Mer14} image \texttt{statiskit/FPD18} (see \url{https://hub.docker.com/r/statiskit/FPD18} for image and the documentation for more information).

\begin{landscape}

\begin{table}
\begin{tabular}{|l|l|l|l|}
\hline
Name 
& Standard beta binomial
& Standard beta negative binomial
& beta Poisson
\\
\hline
Definition
& $\mathcal{B}_{n}\left( p\right) \underset{p}{\wedge} \beta\left(a, b\right)$
& $\mathcal{NB}\left(r, p\right) \underset{p}{\wedge} \beta\left(a, b\right)$
& $\mathcal{P}\left(\lambda p\right) \underset{p}{\wedge} \beta\left(a, b\right)$
\\
Notation
& $\beta\mathcal{B}_{n}\left(a, b\right)$ 
& $\beta\mathcal{NB}\left(r, a, b\right)$
& $\beta\mathcal{P}_{\lambda}\left(a, b\right)$
\\
\hline
Re-parametrization
&  
& UGWD$(r,b,a)$
& 
\\
\hline
Supp($Y$)
& $\{0,1,\ldots,n\}$
& $\mathbb{N}$
& $\mathbb{N}$
\\
\hline
$P(Y=y)$
& $\frac{(b)_{n}}{(a+b)_{n}} \frac{(-n)_y(a)_y}{(-b-n+1)_y} \frac{1}{y!}$
& $\frac{(a)_{b}}{(a+r)_{b}} \frac{(r)_y(b)_y}{(r+a+b)_y}\frac{1}{y!}$ \footnote{$r$ and $b$ have a symmetric role}
& $P(Y=y)=\frac{(a)_y}{(a+b)_y}\frac{\lambda^y}{y!}{}_1F_1(a+y;a+b+y;-\lambda)$
\\
\hline
$\operatorname{E}(Y)$
& $n\frac{a}{a+b}$
& $r\frac{b}{a-1}$ \footnote{define if $a>1$ and undefined otherwise}
& $\lambda\frac{a}{a+b}$
\\
\hline
$\operatorname{V}(Y)$
& $n\frac{ab(a+b+n)}{(a+b)^2(a+b+1)}$
& $r\frac{b(a+r-1)(a+b-1)}{(a-1)^2(a-2)}$ \footnote{define if $a>2$ and undefined otherwise}
& $\lambda\frac{a}{a+b}\left\lbrace 1+\lambda\frac{b}{(a+b)(a+b+1)} \right\rbrace$
\\
\hline
$G_Y(s)$
& $\frac{(b)_{n}}{(a+b)_{n}} \; {}_2F_1\{(-n,a);-b-n+1;s\}$
& $\frac{(a)_{r}}{(a+b)_{r}} \; {}_2F_1\{(r,b);r+a+b;s\}$
& ${}_1F_1\{a;a+b;\lambda(s-1)\}$
\\
\hline
\end{tabular}
\caption{\label{table:standard:beta:compound} Usual characteristics of the standard beta compound binomial, negative binomial and Poisson distributions.}
\end{table}

\begin{table}
\begin{tabular}{|l|l|l|}
\hline
Name 
& Generalized beta binomial
& Generalized beta negative binomial
\\
\hline
Definition
& $\mathcal{B}_{n}\left( p\right) \underset{p}{\wedge} \beta_{\pi}\left(a, b\right)$
& $\mathcal{NB}\left(r, p\right) \underset{p}{\wedge} \beta_{\pi}\left(a, b\right)$
\\
Notation
& $\beta_{\pi}\mathcal{B}_{n}\left(a, b\right)$ 
& $\beta_{\pi}\mathcal{NB}\left(r, a, b\right)$
\\
\hline
Supp($Y$)
& $\{0,1,\ldots,n\}$
& $\mathbb{N}$
\\
\hline
$P(Y=y)$
& $ \frac{(b)_n}{(a+b)_n} \frac{(-n)_y(a)_y}{(-b-n+1)_y} \frac{\pi^y}{y!} {}_2F_1\{(-n+y, a+y); -b-n+1+y; 1-\pi\}$
& $\frac{(a)_r}{(a+b)_r} \frac{(r)_{y}(b)_{y}}{(r+a+b)_y} \frac{\pi^y}{y!} 
    {}_2F_1\{(r+y, b+y); r+a+b+y; 1-\pi \}$
\\
\hline
$\operatorname{E}(Y)$
& $n\pi\frac{a}{a+b}$
& $r\pi\frac{b}{a-1}$ \footnote{define if $a>1$ and undefined otherwise}
\\
\hline
$\operatorname{V}(Y)$
& $n\pi\frac{a}{a+b} \left\lbrace \pi \frac{b(a+b+n)}{(a+b)(a+b+1)} + 1-\pi \right\rbrace$
& $r\pi\frac{b}{a-1}\left\lbrace \pi \frac{(a+r-1)(a+b-1)}{(a-2)(a-1)} + 1-\pi \right\rbrace$ \footnote{define if $a>2$ and undefined otherwise}
\\
\hline
$G_Y(s)$
& $\frac{(b)_{n}}{(a+b)_{n}} \; {}_2F_1\{(-n,a);-b-n+1;1+\pi(s-1)\}$
& $\frac{(a)_{r}}{(a+b)_{r}} \; {}_2F_1\{(r,b);r+a+b;1+\pi(s-1)\}$
\\
\hline
\end{tabular}
\caption{\label{table:generalized:beta:compound} Usual characteristics of the generalized beta binomial and the generalized beta negative binomial distributions.}
\end{table}

\end{landscape}

\end{document}